\documentclass[12pt,thmsa]{article}
\usepackage{amsfonts}
\usepackage{sw20bams}

\input tcilatex
\begin{document}

\author{Steven Finch, Pascal Sebah and Zai-Qiao Bai}
\title{Odd Entries in Pascal's Trinomial Triangle}
\date{February 19, 2008}
\maketitle

\begin{abstract}
The $n^{\text{th}}$ row of Pascal's trinomial triangle gives coefficients of 
$(1+x+x^2)^n$. Let $g(n)$ denote the number of such coefficients that are
odd. We review Moshe's algorithm for evaluating asymptotics of $g(n)$ --
this involves computing the Lyapunov exponent for certain $2\times 2$ random
matrix products -- and then analyze further examples with more terms and
higher powers of $x$.
\end{abstract}

\footnotetext{
Copyright \copyright\ 2008 by Steven R. Finch. All rights reserved.}Before
discussing trinomials, let us recall well-known results for binomials.
Define $f(n)$ to be the number of odd coefficients in $(1+x)^n$. Let $N$
denote a uniform random integer between $0$ and $n-1$, then $f(N)$ has
``typical growth'' $\approx n^{1/2}$ in the sense that 
\[
\limfunc{E}(\ln (f(N)))\sim \frac 12\ln (n) 
\]
as $n\rightarrow \infty ;$ equivalently, 
\[
\lim_{n\rightarrow \infty }\frac 1{n\ln (n)}\dsum\limits_{k=0}^{n-1}\ln
(f(k))=\frac 12=0.5. 
\]
Also $f(N)$ has ``average growth'' $\approx n^{\ln (3/2)/\ln (2)}$ in the
sense that 
\[
\ln (\limfunc{E}(f(N)))\sim \frac{\ln (3/2)}{\ln (2)}\ln (n) 
\]
as $n\rightarrow \infty ;$ equivalently, 
\[
\lim_{n\rightarrow \infty }\frac 1{\ln (n)}\ln \left( \frac
1n\dsum\limits_{k=0}^{n-1}f(k)\right) =\frac{\ln (3/2)}{\ln (2)}%
=0.5849625007211561814537389.... 
\]
The latter value is larger since most of the $1$s in Pascal's binomial
triangle, modulo $2$, are concentrated in relatively few rows. Exact results
are available, due to Trollope \cite{Tr} \& Delange \cite{De} for $\limfunc{E%
}(\ln (f(N)))$ and Stein \cite{St} \&\ Larcher \cite{La} for $\ln (\limfunc{E%
}(f(N)))$; an overview of the subject is found in \cite{Fi}. Our interest
here is solely in first-order approximations.

Moshe \cite{M5} introduced an algorithm for evaluating such asymptotics.
Define $P(i,j)$ to be the $(i,j)^{\text{th}}$ entry of the triangle mod $2$: 
\[
\left( 
\begin{array}{ccccc}
P(0,0) &  &  &  &  \\ 
P(1,0) & P(1,1) &  &  &  \\ 
P(2,0) & P(2,1) & P(2,2) &  &  \\ 
P(3,0) & P(3,1) & P(3,2) & P(3,3) &  \\ 
P(4,0) & P(4,1) & P(4,2) & P(4,3) & P(4,4)
\end{array}
\right) =\left( 
\begin{array}{ccccc}
1 &  &  &  &  \\ 
1 & 1 &  &  &  \\ 
1 & 0 & 1 &  &  \\ 
1 & 1 & 1 & 1 &  \\ 
1 & 0 & 0 & 0 & 1
\end{array}
\right) . 
\]
All values of $P(i,j)$ in the upper right portion are $0$s. Define $o=0$,
which serves as a placeholder. Let $\ell =1$, 
\[
\begin{array}{ccc}
A_0(j)=P(o,j), &  & 0\leq j<\ell
\end{array}
\]
and observe that the vector $A_0=(1)$. Constructing additional $\ell $%
-vectors $A_1$, $A_2$, $\ldots $, $A_{m-1}$, if required, is one aspect of
the algorithm. It is mandatory that $A_i(0)=1$ always. For each $0\leq s\leq
1$, $0\leq t\leq 1$, let 
\[
\begin{array}{ccc}
b_0^{s,t}(j)=P(2o+s,2j+t), &  & 0\leq j<\ell
\end{array}
\]
which is obtained by replacing $(o,j)$ in the expression for $A_0$ by $%
(2o+s,2j+t)$. Observe that the four vectors 
\[
\begin{array}{ccccccc}
b_0^{0,0}=(1), &  & b_0^{1,0}=(1), &  & b_0^{0,1}=(0), &  & b_0^{1,1}=(1)
\end{array}
\]
encompass only $(1)=A_0$ and $(0)$. No ``refinement'' of $b_0^{s,t}$ is
hence necessary (this will be clarified later) and we let 
\[
\begin{array}{ccc}
B_0^{s,t}(j)=P(2o+s,2j+t), &  & 0\leq j<\ell .
\end{array}
\]
Constructing additional $\ell $-vectors $B_1^{s,t}$, $B_2^{s,t}$, $\ldots $, 
$B_{m-1}^{s,t}$, if required, is another aspect of the algorithm. We have 
\[
\begin{array}{ccccccc}
B_0^{0,0}=A_0, &  & B_0^{1,0}=A_0, &  & B_0^{0,1}=0, &  & B_0^{1,1}=A_0
\end{array}
\]
and this completes the iterative portion of the algorithm. Thus $m=1$.
Define $m\times m$ matrices via 
\[
D_s(i,j)=\#\left\{ t:B_j^{s,t}=A_i\right\} 
\]
and observe that $D_0=(1)$, $D_1=(2)$. Define $m$-vectors via 
\[
\begin{array}{ccc}
e_i(j)=\left\{ 
\begin{array}{lll}
1 &  & \text{if }j=i, \\ 
0 &  & \text{if }j\neq i,
\end{array}
\right. &  & w(i)=\#\left\{ j:A_i(j)=1\right\}
\end{array}
\]
and observe that $e_0=(1)$, $w=(1)$. Let the binary expansion of a positive
integer $n$ be $\sum_{r=0}^{k-1}n_r2^r$ with $0\leq n_r<2$ and $n_{k-1}=1$.
Moshe \cite{M5} proved the following formula: 
\[
f(n)=\#\left\{ j:P(n,j)=1\right\} =w^TD_{n_{k-1}}D_{n_{k-2}}\cdots
D_{n_1}D_{n_0}e_0 
\]
which provides a useful check that $D_0$, $D_1$, $w$ are correct.

One consequence is a well-known recursive formula for $f(n)$. Writing binary
expansions as $n=n_{k-1}n_{k-2}\ldots n_1n_0$, we see 
\[
\begin{array}{ccc}
2n=n_{k-1}n_{k-2}\ldots n_1n_00, &  & 2n+1=n_{k-1}n_{k-2}\ldots n_1n_01
\end{array}
\]
hence 
\begin{eqnarray*}
f(2n) &=&w^TD_{n_{k-1}}D_{n_{k-2}}\cdots D_{n_1}D_{n_0}D_0e_0 \\
&=&w^TD_{n_{k-1}}D_{n_{k-2}}\cdots D_{n_1}D_{n_0}e_0 \\
&=&f(n),
\end{eqnarray*}
\begin{eqnarray*}
f(2n+1) &=&w^TD_{n_{k-1}}D_{n_{k-2}}\cdots D_{n_1}D_{n_0}D_1e_0 \\
&=&2w^TD_{n_{k-1}}D_{n_{k-2}}\cdots D_{n_1}D_{n_0}e_0 \\
&=&2f(n).
\end{eqnarray*}
All $D_0$ matrices exhibited in this paper satisfy $D_0e_0=e_0$, therefore
even arguments are easy. Odd arguments are harder since $D_1e_0$ is not as
predictable.

Another consequence involves the growth rates $\limfunc{E}(\ln (f(N)))$ and $%
\ln (\limfunc{E}(f(N)))$. Let us work with the latter first. Order the
complex eigenvalues $\mu _1$, $\mu _2$, $\ldots $, $\mu _m$ of the matrix $%
D_0+D_1$ so that $\mu _1$ has maximum modulus; therefore \cite{M5} 
\[
\frac{\ln (\limfunc{E}(f(N)))}{\ln (n)}\rightarrow \frac{\ln \left| \mu
_1\right| }{\ln (2)}-1=\frac{\ln (3/2)}{\ln (2)} 
\]
as $n\rightarrow \infty $. Working with $\limfunc{E}(\ln (f(N)))$ is more
complicated. In this scalar case, it is clear that for $2^{k-1}\leq N<n\leq
2^k$, 
\[
\limfunc{E}(\ln (f(N)))=\limfunc{E}\left( \dsum_{r=0}^{k-1}\ln
(D_{N_r})\right) =\left( \frac 12\ln (1)+\frac 12\ln (2)\right) k=\frac{\ln
(2)}2\left\lceil \frac{\ln (n)}{\ln (2)}\right\rceil 
\]
and so 
\[
\frac{\limfunc{E}(\ln (f(N)))}{\ln (n)}\rightarrow \frac 12 
\]
as $n\rightarrow \infty $. But commutativity fails for random matrix
products, in general, and a Lyapunov exponent-based approach will be
presented in section [\ref{Lyap}] to deal with this issue.

\subsection{Trinomials I}

Define $g(n)$ to be the number of odd coefficients in $(1+x+x^2)^n$.
Properties of $g(n)$ and Pascal's trinomial triangle are given in \cite{Eu,
Bo}. Define $P(i,j)$ to be the $(i,j)^{\text{th}}$ entry of the triangle mod 
$2$: 
\[
\left( 
\begin{array}{cccccccc}
\QTR{textit}{{\scriptsize {P}}}\text{{\scriptsize (0,0)}} & \QTR{textit}{%
{\scriptsize {P}}}\text{{\scriptsize (0,1)}} &  &  &  &  &  &  \\ 
\QTR{textit}{{\scriptsize {P}}}\text{{\scriptsize (1,0)}} & \QTR{textit}{%
{\scriptsize {P}}}\text{{\scriptsize (1,1)}} & \QTR{textit}{{\scriptsize {P}}%
}\text{{\scriptsize (1,2)}} & \QTR{textit}{{\scriptsize {P}}}\text{%
{\scriptsize (1,3)}} &  &  &  &  \\ 
\QTR{textit}{{\scriptsize {P}}}\text{{\scriptsize (2,0)}} & \QTR{textit}{%
{\scriptsize {P}}}\text{{\scriptsize (2,1)}} & \QTR{textit}{{\scriptsize {P}}%
}\text{{\scriptsize (2,2)}} & \QTR{textit}{{\scriptsize {P}}}\text{%
{\scriptsize (2,3)}} & \QTR{textit}{{\scriptsize {P}}}\text{{\scriptsize %
(2,4)}} & \QTR{textit}{{\scriptsize {P}}}\text{{\scriptsize (2,5)}} &  &  \\ 
\QTR{textit}{{\scriptsize {P}}}\text{{\scriptsize (3,0)}} & \QTR{textit}{%
{\scriptsize {P}}}\text{{\scriptsize (3,1)}} & \QTR{textit}{{\scriptsize {P}}%
}\text{{\scriptsize (3,2)}} & \QTR{textit}{{\scriptsize {P}}}\text{%
{\scriptsize (3,3)}} & \QTR{textit}{{\scriptsize {P}}}\text{{\scriptsize %
(3,4)}} & \QTR{textit}{{\scriptsize {P}}}\text{{\scriptsize (3,5)}} & 
\QTR{textit}{{\scriptsize {P}}}\text{{\scriptsize (3,6)}} & \QTR{textit}{%
{\scriptsize {P}}}\text{{\scriptsize (3,7)}}
\end{array}
\right) =\left( 
\begin{array}{cccccccc}
\text{{\scriptsize 1}} & \text{{\scriptsize 0}} &  &  &  &  &  &  \\ 
\text{{\scriptsize 1}} & \text{{\scriptsize 1}} & \text{{\scriptsize 1}} & 
\text{{\scriptsize 0}} &  &  &  &  \\ 
\text{{\scriptsize 1}} & \text{{\scriptsize 0}} & \text{{\scriptsize 1}} & 
\text{{\scriptsize 0}} & \text{{\scriptsize 1}} & \text{{\scriptsize 0}} & 
&  \\ 
\text{{\scriptsize 1}} & \text{{\scriptsize 1}} & \text{{\scriptsize 0}} & 
\text{{\scriptsize 1}} & \text{{\scriptsize 0}} & \text{{\scriptsize 1}} & 
\text{{\scriptsize 1}} & \text{{\scriptsize 0}}
\end{array}
\right) . 
\]
All values of $P(i,j)$ in the upper right portion are $0$s. Let $\ell =2$, 
\[
\begin{array}{ccc}
A_0(j)=P(o,j), &  & 0\leq j<\ell
\end{array}
\]
and observe that the vector $A_0=(1,0)$. For each $0\leq s\leq 1$, $0\leq
t\leq 1$, let 
\[
\begin{array}{ccc}
b_0^{s,t}(j)=P(2o+s,2j+t), &  & 0\leq j<\ell
\end{array}
\]
which is obtained by replacing $(o,j)$ in the expression for $A_0$ by $%
(2o+s,2j+t)$. Observe that the four vectors 
\[
\begin{array}{ccccccc}
b_0^{0,0}=(1,0), &  & b_0^{1,0}=(1,1), &  & b_0^{0,1}=(0,0), &  & 
b_0^{1,1}=(1,0)
\end{array}
\]
contain a nonzero vector $(1,1)\neq A_0$. Let $A_1=b_0^{1,0}=(1,1)$. No
``refinement'' of $b_0^{s,t}$ is necessary (this will be clarified soon) and
we let 
\[
\begin{array}{ccc}
B_0^{s,t}(j)=P(2o+s,2j+t), &  & 0\leq j<\ell .
\end{array}
\]
It follows that 
\[
\begin{array}{ccccccc}
B_0^{0,0}=A_0, &  & B_0^{1,0}=A_1, &  & B_0^{0,1}=0, &  & B_0^{1,1}=A_0
\end{array}
\]
but we are not yet done (because of $A_1$). Let 
\[
\begin{array}{ccc}
b_1^{s,t}(j)=P(2(2o+s)+1,2(2j+t)+0), &  & 0\leq j<\ell
\end{array}
\]
which is obtained by replacing $(o,j)$ in the expression for $b_0^{1,0}=A_1$
by $(2o+s,2j+t)$. Observe that the four vectors 
\[
\begin{array}{ccccccc}
b_1^{0,0}=(1,0), &  & b_1^{1,0}=(1,0), &  & b_1^{0,1}=(1,0), &  & 
b_1^{1,1}=(0,1)
\end{array}
\]
contain a nonzero vector $(0,1)\neq A_0,A_1$. Setting $A_2=b_1^{1,1}$,
however, violates the mandate that $A_i(0)=1$ always. We thus refine the
definition of $b_1^{s,t}$: 
\[
B_1^{s,t}(j)=\left\{ 
\begin{array}{lll}
P(2(2o+s)+1,2(2j+t)+0) &  & \text{if }(s,t)\neq (1,1), \\ 
P(2(2o+s)+1,2(2(j+c)+t)+0) &  & \text{if }(s,t)=(1,1)
\end{array}
\right. 
\]
where $c=1$, which corresponds to shifting one step to the right. It follows
that 
\[
\begin{array}{ccccccc}
B_1^{0,0}=A_0, &  & B_1^{1,0}=A_0, &  & B_1^{0,1}=A_0, &  & B_1^{1,1}=A_0
\end{array}
\]
and this completes the iterative portion of the algorithm. Thus $m=2$. The
definitions of $m\times m$ matrices 
\[
\begin{array}{ccc}
D_0=\left( 
\begin{array}{cc}
1 & 2 \\ 
0 & 0
\end{array}
\right) , &  & D_1=\left( 
\begin{array}{cc}
1 & 2 \\ 
1 & 0
\end{array}
\right) .
\end{array}
\]
and $m$-vectors $e_0=(1,0)$, $w=(1,2)$ are as before. Likewise, 
\[
g(n)=\#\left\{ j:P(n,j)=1\right\} =w^TD_{n_{k-1}}D_{n_{k-2}}\cdots
D_{n_1}D_{n_0}e_0 
\]
as before, where $n_{k-1}n_{k-2}\ldots n_1n_0$ is the binary expansion of $n$%
.

It is easy to see that $g(2n)=g(n)$. From 
\[
\begin{array}{ccc}
2n+1=n_{k-1}n_{k-2}\ldots n_1n_01, &  & D_1e_0=e_0+e_1,
\end{array}
\]
\[
\begin{array}{ccc}
4n+1=n_{k-1}n_{k-2}\ldots n_1n_001, &  & D_0D_1e_0=3e_0,
\end{array}
\]
\[
\begin{array}{ccc}
4n+3=n_{k-1}n_{k-2}\ldots n_1n_011, &  & D_1D_1e_0=3e_0+e_1
\end{array}
\]
we reproduce Sillke's result \cite{Si} that 
\[
\begin{array}{ccc}
g(4n+1)=3g(n), &  & g(4n+3)=g(2n+1)+2g(n).
\end{array}
\]
Also, the maximal eigenvalue $\mu _1$ of $D_0+D_1$ is $1+\sqrt{5}$;
therefore 
\[
\frac{\ln (\limfunc{E}(g(N)))}{\ln (n)}\rightarrow \frac{\ln (1+\sqrt{5})}{%
\ln (2)}-1=\frac{\ln (\varphi )}{\ln (2)}=0.6942419136306173017387902... 
\]
as $n\rightarrow \infty $, where $\varphi $ is the Golden mean. This
constant is not new: see \cite{Wf}.

\subsection{Quadrinomials}

Define $g_3(n)$ to be the number of odd coefficients in $(1+x+x^2+x^3)^n$.
This extends our earlier definitions $f=g_1$ and $g=g_2$. Define $P(i,j)$ to
be the $(i,j)^{\text{th}}$ entry of Pascal's quadrinomial triangle mod $2$: 
\[
\left( 
\begin{array}{cccccc}
P(0,0) & P(0,1) & P(0,2) &  &  &  \\ 
P(1,0) & P(1,1) & P(1,2) & P(1,3) & P(1,4) & P(1,5)
\end{array}
\right) =\left( 
\begin{array}{cccccc}
1 & 0 & 0 &  &  &  \\ 
1 & 1 & 1 & 1 & 0 & 0
\end{array}
\right) . 
\]
All values of $P(i,j)$ in the upper right portion are $0$s. Let $\ell =3$, 
\[
\begin{array}{ccc}
A_0(j)=P(o,j), &  & 0\leq j<\ell
\end{array}
\]
and observe that $A_0=(1,0,0)$. For each $0\leq s\leq 1$, $0\leq t\leq 1$,
let 
\[
\begin{array}{ccc}
b_0^{s,t}(j)=P(2o+s,2j+t), &  & 0\leq j<\ell
\end{array}
\]
and observe that the four vectors 
\[
\begin{array}{ccccccc}
b_0^{0,0}=(1,0,0), &  & b_0^{1,0}=(1,1,0), &  & b_0^{0,1}=(0,0,0), &  & 
b_0^{1,1}=(1,1,0)
\end{array}
\]
contain a nonzero vector $(1,1,0)\neq A_0$. Let $A_1=b_0^{1,0}=(1,1,0)$. No
refinement of $b_0^{s,t}$ is necessary and we let 
\[
\begin{array}{ccc}
B_0^{s,t}(j)=P(2o+s,2j+t), &  & 0\leq j<\ell .
\end{array}
\]
It follows that 
\[
\begin{array}{ccccccc}
B_0^{0,0}=A_0, &  & B_0^{1,0}=A_1, &  & B_0^{0,1}=0, &  & B_0^{1,1}=A_1
\end{array}
\]
but we are not yet done. Let 
\[
\begin{array}{ccc}
b_1^{s,t}(j)=P(2(2o+s)+1,2(2j+t)+0), &  & 0\leq j<\ell
\end{array}
\]
which is obtained by replacing $(o,j)$ in the expression for $b_0^{1,0}=A_1$
by $(2o+s,2j+t)$. Observe that 
\[
\begin{array}{ccccccc}
b_1^{0,0}=(1,0,0), &  & b_1^{1,0}=(1,0,1), &  & b_1^{0,1}=(1,0,0), &  & 
b_1^{1,1}=(0,0,0)
\end{array}
\]
contain a nonzero vector $(1,0,1)\neq A_0,A_1$. Let $A_2=b_1^{1,0}=(1,0,1)$.
No refinement of $b_1^{s,t}$ is necessary and we let 
\[
\begin{array}{ccc}
B_1^{s,t}(j)=P(2(2o+s)+1,2(2j+t)+0), &  & 0\leq j<\ell .
\end{array}
\]
It follows that 
\[
\begin{array}{ccccccc}
B_1^{0,0}=A_0, &  & B_1^{1,0}=A_2, &  & B_1^{0,1}=A_0, &  & B_1^{1,1}=0
\end{array}
\]
but we are not yet done. Let 
\[
\begin{array}{ccc}
b_2^{s,t}(j)=P(2(2(2o+s)+1)+1,2(2(2j+t)+0)+0), &  & 0\leq j<\ell
\end{array}
\]
which is obtained by replacing $(o,j)$ in the expression for $b_1^{1,0}=A_2$
by $(2o+s,2j+t)$. Observe that 
\[
\begin{array}{ccccccc}
b_2^{0,0}=(1,1,0), &  & b_2^{1,0}=(1,0,1), &  & b_2^{0,1}=(0,0,0), &  & 
b_2^{1,1}=(1,0,1)
\end{array}
\]
encompass only $A_1$, $A_2$ and $0$. No refinement of $b_2^{s,t}$ is hence
necessary and we let 
\[
\begin{array}{ccc}
B_2^{s,t}(j)=P(2(2(2o+s)+1)+1,2(2(2j+t)+0)+0), &  & 0\leq j<\ell .
\end{array}
\]
It follows that 
\[
\begin{array}{ccccccc}
B_2^{0,0}=A_1, &  & B_2^{1,0}=A_2, &  & B_2^{0,1}=0, &  & B_2^{1,1}=A_2
\end{array}
\]
and this completes the iterative portion of the algorithm. Thus $m=3$. The
definitions of $m\times m$ matrices 
\[
\begin{array}{ccc}
D_0=\left( 
\begin{array}{ccc}
1 & 2 & 0 \\ 
0 & 0 & 1 \\ 
0 & 0 & 0
\end{array}
\right) , &  & D_1=\left( 
\begin{array}{ccc}
0 & 0 & 0 \\ 
2 & 0 & 0 \\ 
0 & 1 & 2
\end{array}
\right)
\end{array}
\]
and $m$-vectors $e_0=(1,0,0)$, $w=(1,2,2)$ are as before. Likewise, 
\[
g_3(n)=\#\left\{ j:P(n,j)=1\right\} =w^TD_{n_{k-1}}D_{n_{k-2}}\cdots
D_{n_1}D_{n_0}e_0 
\]
as before, where $n_{k-1}n_{k-2}\ldots n_1n_0$ is the binary expansion of $n$%
.

It is easy to see that $g_3(2n)=g_3(n)$. From 
\[
\begin{array}{ccc}
2n+1=n_{k-1}n_{k-2}\ldots n_1n_01, &  & D_1e_0=2e_1,
\end{array}
\]
\[
\begin{array}{ccc}
4n+1=n_{k-1}n_{k-2}\ldots n_1n_001, &  & D_0D_1e_0=4e_0,
\end{array}
\]
\[
\begin{array}{ccc}
4n+3=n_{k-1}n_{k-2}\ldots n_1n_011, &  & D_1D_1e_0=2e_2,
\end{array}
\]
\[
\begin{array}{ccc}
8n+1=n_{k-1}n_{k-2}\ldots n_1n_0001, &  & D_0D_0D_1e_0=4e_0,
\end{array}
\]
\[
\begin{array}{ccc}
8n+3=n_{k-1}n_{k-2}\ldots n_1n_0011, &  & D_0D_1D_1e_0=2e_1,
\end{array}
\]
\[
\begin{array}{ccc}
8n+5=n_{k-1}n_{k-2}\ldots n_1n_0101, &  & D_1D_0D_1e_0=8e_1,
\end{array}
\]
\[
\begin{array}{ccc}
8n+7=n_{k-1}n_{k-2}\ldots n_1n_0111, &  & D_1D_1D_1e_0=4e_2
\end{array}
\]
we deduce that 
\[
\begin{array}{ccc}
g_3(8n+1)=g_3(4n+1), &  & g_3(8n+3)=g_3(2n+1),
\end{array}
\]
\[
\begin{array}{ccc}
g_3(8n+5)=4g_3(2n+1), &  & g_3(8n+7)=2g_3(4n+3).
\end{array}
\]
Also, the maximal eigenvalue $\mu _1$ of $D_0+D_1$ is $3$; it is interesting
that the same constant 
\[
\frac{\ln (\limfunc{E}(g_3(N)))}{\ln (n)}\rightarrow \frac{\ln (3/2)}{\ln (2)%
} 
\]
appears here as for $g_1(N).$

\subsection{Trinomials II}

Define $h_3(n)$ to be the number of odd coefficients in $(1+x+x^3)^n$. This
extends our earlier definition $g=h_2$. Define $P(i,j)$ to be the $(i,j)^{%
\text{th}}$ entry of the associated triangle mod $2$: 
\[
\left( 
\begin{array}{cccccc}
P(0,0) & P(0,1) & P(0,2) &  &  &  \\ 
P(1,0) & P(1,1) & P(1,2) & P(1,3) & P(1,4) & P(1,5)
\end{array}
\right) =\left( 
\begin{array}{cccccc}
1 & 0 & 0 &  &  &  \\ 
1 & 1 & 0 & 1 & 0 & 0
\end{array}
\right) . 
\]
All values of $P(i,j)$ in the upper right portion are $0$s. Let $\ell =3$, 
\[
\begin{array}{ccc}
A_0(j)=P(o,j), &  & 0\leq j<\ell
\end{array}
\]
and observe that $A_0=(1,0,0)$. For each $0\leq s\leq 1$, $0\leq t\leq 1$,
let 
\[
\begin{array}{ccc}
b_0^{s,t}(j)=P(2o+s,2j+t), &  & 0\leq j<\ell
\end{array}
\]
and observe that the four vectors 
\[
\begin{array}{ccccccc}
b_0^{0,0}=(1,0,0), &  & b_0^{1,0}=(1,0,0), &  & b_0^{0,1}=(0,0,0), &  & 
b_0^{1,1}=(1,1,0)
\end{array}
\]
contain a nonzero vector $(1,1,0)\neq A_0$. Let $A_1=b_0^{1,1}=(1,1,0)$. No
refinement of $b_0^{s,t}$ is necessary and we let 
\[
\begin{array}{ccc}
B_0^{s,t}(j)=P(2o+s,2j+t), &  & 0\leq j<\ell .
\end{array}
\]
It follows that 
\[
\begin{array}{ccccccc}
B_0^{0,0}=A_0, &  & B_0^{1,0}=A_0, &  & B_0^{0,1}=0, &  & B_0^{1,1}=A_1
\end{array}
\]
but we are not yet done. Let 
\[
\begin{array}{ccc}
b_1^{s,t}(j)=P(2(2o+s)+1,2(2j+t)+1), &  & 0\leq j<\ell
\end{array}
\]
which is obtained by replacing $(o,j)$ in the expression for $b_0^{1,1}=A_1$
by $(2o+s,2j+t)$. Observe that 
\[
\begin{array}{ccccccc}
b_1^{0,0}=(1,0,0), &  & b_1^{1,0}=(1,1,1), &  & b_1^{0,1}=(1,0,0), &  & 
b_1^{1,1}=(0,1,0)
\end{array}
\]
contain nonzero vectors $(1,1,1),(0,1,0)\neq A_0,A_1$. Let $%
A_2=b_1^{1,0}=(1,1,1)$. Setting $A_3=b_1^{1,1}$, however, violates the
mandate that $A_i(0)=1$ always. We thus refine the definition of $b_1^{s,t}$%
: 
\[
B_1^{s,t}(j)=\left\{ 
\begin{array}{lll}
P(2(2o+s)+1,2(2j+t)+1) &  & \text{if }(s,t)\neq (1,1), \\ 
P(2(2o+s)+1,2(2(j+c)+t)+1) &  & \text{if }(s,t)=(1,1)
\end{array}
\right. 
\]
where $c=1$. It follows that 
\[
\begin{array}{ccccccc}
B_1^{0,0}=A_0, &  & B_1^{1,0}=A_2, &  & B_1^{0,1}=A_0, &  & B_1^{1,1}=A_0
\end{array}
\]
but we are not yet done. Let 
\[
\begin{array}{ccc}
b_2^{s,t}(j)=P(2(2(2o+s)+1)+1,2(2(2j+t)+0)+1), &  & 0\leq j<\ell
\end{array}
\]
which is obtained by replacing $(o,j)$ in the expression for $b_1^{1,0}=A_2$
by $(2o+s,2j+t)$. Observe that 
\[
\begin{array}{ccccccc}
b_2^{0,0}=(1,1,0), &  & b_2^{1,0}=(1,0,1), &  & b_2^{0,1}=(1,0,0), &  & 
b_2^{1,1}=(0,0,1)
\end{array}
\]
contain nonzero vectors $(1,0,1),(0,0,1)\neq A_0,A_1,A_2$. Let $%
A_3=b_2^{1,0}=(1,0,1)$. Setting $A_4=b_2^{1,1}$, however, violates the
mandate that $A_i(0)=1$ always. We thus refine the definition of $b_2^{s,t}$%
: 
\[
B_2^{s,t}(j)=\left\{ 
\begin{array}{lll}
P(2(2(2o+s)+1)+1,2(2(2j+t)+0)+1) &  & \text{if }(s,t)\neq (1,1), \\ 
P(2(2(2o+s)+1)+1,2(2(2(j+c)+t)+0)+1) &  & \text{if }(s,t)=(1,1)
\end{array}
\right. 
\]
where $c=2$. It follows that 
\[
\begin{array}{ccccccc}
B_2^{0,0}=A_1, &  & B_2^{1,0}=A_3, &  & B_2^{0,1}=A_0, &  & B_2^{1,1}=A_0
\end{array}
\]
but we are not yet done. Let 
\[
\begin{array}{ccc}
b_3^{s,t}(j)=P(2(2(2(2o+s)+1)+1)+1,2(2(2(2j+t)+0)+0)+1), &  & 0\leq j<\ell
\end{array}
\]
which is obtained by replacing $(o,j)$ in the expression for $b_2^{1,0}=A_3$
by $(2o+s,2j+t)$. Observe that 
\[
\begin{array}{ccccccc}
b_3^{0,0}=(1,1,0), &  & b_3^{1,0}=(1,1,0), &  & b_3^{0,1}=(0,0,0), &  & 
b_3^{1,1}=(1,0,1)
\end{array}
\]
encompass only $A_1$, $A_3$ and $0$. No refinement of $b_3^{s,t}$ is hence
necessary and we let 
\[
\begin{array}{ccc}
B_3^{s,t}(j)=P(2(2(2(2o+s)+1)+1)+1,2(2(2(2j+t)+0)+0)+1), &  & 0\leq j<\ell .
\end{array}
\]
It follows that 
\[
\begin{array}{ccccccc}
B_3^{0,0}=A_1, &  & B_3^{1,0}=A_1, &  & B_3^{0,1}=0, &  & B_3^{1,1}=A_3
\end{array}
\]
and this completes the iterative portion of the algorithm. Thus $m=4$. The
definitions of $m\times m$ matrices 
\[
\begin{array}{ccc}
D_0=\left( 
\begin{array}{cccc}
1 & 2 & 1 & 0 \\ 
0 & 0 & 1 & 1 \\ 
0 & 0 & 0 & 0 \\ 
0 & 0 & 0 & 0
\end{array}
\right) , &  & D_1=\left( 
\begin{array}{cccc}
1 & 1 & 1 & 0 \\ 
1 & 0 & 0 & 1 \\ 
0 & 1 & 0 & 0 \\ 
0 & 0 & 1 & 1
\end{array}
\right)
\end{array}
\]
and $m$-vectors $e_0=(1,0,0,0)$, $w=(1,2,3,2)$ are as before. Likewise, 
\[
h_3(n)=\#\left\{ j:P(n,j)=1\right\} =w^TD_{n_{k-1}}D_{n_{k-2}}\cdots
D_{n_1}D_{n_0}e_0 
\]
as before, where $n_{k-1}n_{k-2}\ldots n_1n_0$ is the binary expansion of $n$%
.

It is easy to see that $h_3(2n)=h_3(n)$. Omitting details, we deduce that 
\[
\begin{array}{ccc}
h_3(4n+1)=3h_3(n), &  & h_3(8n+3)=h_3(2n+1)+4h_3(n),
\end{array}
\]
\[
h_3(16n+7)=h_3(8n+3)+h_3(2n+1)+\,3h_3(n), 
\]
\[
h_3(16n+15)=2h_3(8n+7)+h_3(2n+1)-2h_3(n). 
\]
Also, the eigenvalues of $D_0+D_1$ have minimal polynomial $\xi ^4-3\xi
^3-2\xi ^2+2\xi +4$ and 
\[
\frac{\ln (\limfunc{E}(h_3(N)))}{\ln (n)}\rightarrow
0.7274509132400228143266172... 
\]
as $n\rightarrow \infty $.

\subsection{Quintinomials}

Define $g_4(n)$ to be the number of odd coefficients in $(1+x+x^2+x^3+x^4)^n$%
. Define $P(i,j)$ to be the $(i,j)^{\text{th}}$ entry of Pascal's
quintinomial triangle mod $2$: 
\[
\left( 
\begin{array}{cccccccc}
\QTR{textit}{{\scriptsize {P}}}\text{{\scriptsize (0,0)}} & \QTR{textit}{%
{\scriptsize {P}}}\text{{\scriptsize (0,1)}} & \QTR{textit}{{\scriptsize {P}}%
}\text{{\scriptsize (0,2)}} & \QTR{textit}{{\scriptsize {P}}}\text{%
{\scriptsize (0,3)}} &  &  &  &  \\ 
\QTR{textit}{{\scriptsize {P}}}\text{{\scriptsize (1,0)}} & \QTR{textit}{%
{\scriptsize {P}}}\text{{\scriptsize (1,1)}} & \QTR{textit}{{\scriptsize {P}}%
}\text{{\scriptsize (1,2)}} & \QTR{textit}{{\scriptsize {P}}}\text{%
{\scriptsize (1,3)}} & \QTR{textit}{{\scriptsize {P}}}\text{{\scriptsize %
(1,4)}} & \QTR{textit}{{\scriptsize {P}}}\text{{\scriptsize (1,5)}} & 
\QTR{textit}{{\scriptsize {P}}}\text{{\scriptsize (1,6)}} & \QTR{textit}{%
{\scriptsize {P}}}\text{{\scriptsize (1,7)}}
\end{array}
\right) =\left( 
\begin{array}{cccccccc}
\text{{\scriptsize 1}} & \text{{\scriptsize 0}} & \text{{\scriptsize 0}} & 
\text{{\scriptsize 0}} &  &  &  &  \\ 
\text{{\scriptsize 1}} & \text{{\scriptsize 1}} & \text{{\scriptsize 1}} & 
\text{{\scriptsize 1}} & \text{{\scriptsize 1}} & \text{{\scriptsize 0}} & 
\text{{\scriptsize 0}} & \text{{\scriptsize 0}}
\end{array}
\right) . 
\]
Let $\ell =4$ and $A_0=(1,0,0,0)$. We simply summarize: 
\[
\begin{array}{ccc}
A_1=b_0^{1,0}=(1,1,1,0), &  & A_2=b_0^{1,1}=(1,1,0,0);
\end{array}
\]
\[
\begin{array}{ccccccc}
B_0^{0,0}=A_0, &  & B_0^{1,0}=A_1, &  & B_0^{0,1}=0, &  & B_0^{1,1}=A_2;
\end{array}
\]
\[
A_3=b_1^{1,0}=(1,1,1,1); 
\]
\[
\begin{array}{ccccccc}
B_1^{0,0}=A_2, &  & B_1^{1,0}=A_3, &  & B_1^{0,1}=A_0, &  & B_1^{1,1}=A_0;
\end{array}
\]
\[
\begin{array}{ccccccc}
B_2^{0,0}=A_0, &  & B_2^{1,0}=A_0, &  & B_2^{0,1}=A_0, &  & B_2^{1,1}=A_0;
\end{array}
\]
\[
\begin{array}{ccccccc}
B_3^{0,0}=A_2, &  & B_3^{1,0}=A_2, &  & B_3^{0,1}=A_2, &  & B_3^{1,1}=A_2;
\end{array}
\]
\[
\begin{array}{ccc}
D_0=\left( 
\begin{array}{cccc}
1 & 1 & 2 & 0 \\ 
0 & 0 & 0 & 0 \\ 
0 & 1 & 0 & 2 \\ 
0 & 0 & 0 & 0
\end{array}
\right) , &  & D_1=\left( 
\begin{array}{cccc}
0 & 1 & 2 & 0 \\ 
1 & 0 & 0 & 0 \\ 
1 & 0 & 0 & 2 \\ 
0 & 1 & 0 & 0
\end{array}
\right)
\end{array}
\]
and $w=(1,3,2,4)$.

The eigenvalues of $D_0+D_1$ have minimal polynomial $\xi ^4-\xi ^3-6\xi
^2-4\xi -16$ and 
\[
\frac{\ln (\limfunc{E}(g_4(N)))}{\ln (n)}\rightarrow
0.7896418505307685639015472... 
\]
as $n\rightarrow \infty $.

\subsection{Trinomials III}

Define $h_4(n)$ to be the number of odd coefficients in $(1+x+x^4)^n$.
Define $P(i,j)$ to be the $(i,j)^{\text{th}}$ entry of the associated
triangle mod $2$: 
\[
\left( 
\begin{array}{cccccccc}
\QTR{textit}{{\scriptsize {P}}}\text{{\scriptsize (0,0)}} & \QTR{textit}{%
{\scriptsize {P}}}\text{{\scriptsize (0,1)}} & \QTR{textit}{{\scriptsize {P}}%
}\text{{\scriptsize (0,2)}} & \QTR{textit}{{\scriptsize {P}}}\text{%
{\scriptsize (0,3)}} &  &  &  &  \\ 
\QTR{textit}{{\scriptsize {P}}}\text{{\scriptsize (1,0)}} & \QTR{textit}{%
{\scriptsize {P}}}\text{{\scriptsize (1,1)}} & \QTR{textit}{{\scriptsize {P}}%
}\text{{\scriptsize (1,2)}} & \QTR{textit}{{\scriptsize {P}}}\text{%
{\scriptsize (1,3)}} & \QTR{textit}{{\scriptsize {P}}}\text{{\scriptsize %
(1,4)}} & \QTR{textit}{{\scriptsize {P}}}\text{{\scriptsize (1,5)}} & 
\QTR{textit}{{\scriptsize {P}}}\text{{\scriptsize (1,6)}} & \QTR{textit}{%
{\scriptsize {P}}}\text{{\scriptsize (1,7)}}
\end{array}
\right) =\left( 
\begin{array}{cccccccc}
\text{{\scriptsize 1}} & \text{{\scriptsize 0}} & \text{{\scriptsize 0}} & 
\text{{\scriptsize 0}} &  &  &  &  \\ 
\text{{\scriptsize 1}} & \text{{\scriptsize 1}} & \text{{\scriptsize 0}} & 
\text{{\scriptsize 0}} & \text{{\scriptsize 1}} & \text{{\scriptsize 0}} & 
\text{{\scriptsize 0}} & \text{{\scriptsize 0}}
\end{array}
\right) . 
\]
Let $\ell =4$ and $A_0=(1,0,0,0)$. We simply summarize: 
\[
A_1=b_0^{1,0}=(1,0,1,0); 
\]
\[
\begin{array}{ccccccc}
B_0^{0,0}=A_0, &  & B_0^{1,0}=A_1, &  & B_0^{0,1}=0, &  & B_0^{1,1}=A_0;
\end{array}
\]
\[
\begin{array}{ccc}
A_2=b_1^{0,0}=(1,1,0,0), &  & A_3=b_1^{1,0}=(1,1,1,1);
\end{array}
\]
\[
\begin{array}{ccccccc}
B_1^{0,0}=A_2, &  & B_1^{1,0}=A_3, &  & B_1^{0,1}=0, &  & B_1^{1,1}=A_2;
\end{array}
\]
\[
A_4=b_2^{1,0}=(1,1,1,0); 
\]
\[
\begin{array}{ccccccc}
B_2^{0,0}=A_0, &  & B_2^{1,0}=A_4, &  & B_2^{0,1}=A_0, &  & B_2^{1,1}=A_0;
\end{array}
\]
\[
A_5=b_3^{1,0}=(1,0,0,1); 
\]
\[
\begin{array}{ccccccc}
B_3^{0,0}=A_2, &  & B_3^{1,0}=A_5, &  & B_3^{0,1}=A_2, &  & B_3^{1,1}=A_2;
\end{array}
\]
\[
A_6=b_4^{1,0}=(1,0,1,1); 
\]
\[
\begin{array}{ccccccc}
B_4^{0,0}=A_2, &  & B_4^{1,0}=A_6, &  & B_4^{0,1}=A_0, &  & B_4^{1,1}=A_2;
\end{array}
\]
\[
A_7=b_5^{1,1}=(1,1,0,1); 
\]
\[
\begin{array}{ccccccc}
B_5^{0,0}=A_0, &  & B_5^{1,0}=A_0, &  & B_5^{0,1}=A_0, &  & B_5^{1,1}=A_7;
\end{array}
\]
\[
\begin{array}{ccccccc}
B_6^{0,0}=A_2, &  & B_6^{1,0}=A_7, &  & B_6^{0,1}=A_0, &  & B_6^{1,1}=A_5;
\end{array}
\]
\[
\begin{array}{ccccccc}
B_7^{0,0}=A_0, &  & B_7^{1,0}=A_2, &  & B_7^{0,1}=A_2, &  & B_7^{1,1}=A_4;
\end{array}
\]
\[
\begin{array}{ccc}
D_0=\left( 
\begin{array}{cccccccc}
\text{{\scriptsize 1}} & \text{{\scriptsize 0}} & \text{{\scriptsize 2}} & 
\text{{\scriptsize 0}} & \text{{\scriptsize 1}} & \text{{\scriptsize 2}} & 
\text{{\scriptsize 1}} & \text{{\scriptsize 1}} \\ 
\text{{\scriptsize 0}} & \text{{\scriptsize 0}} & \text{{\scriptsize 0}} & 
\text{{\scriptsize 0}} & \text{{\scriptsize 0}} & \text{{\scriptsize 0}} & 
\text{{\scriptsize 0}} & \text{{\scriptsize 0}} \\ 
\text{{\scriptsize 0}} & \text{{\scriptsize 1}} & \text{{\scriptsize 0}} & 
\text{{\scriptsize 2}} & \text{{\scriptsize 1}} & \text{{\scriptsize 0}} & 
\text{{\scriptsize 1}} & \text{{\scriptsize 1}} \\ 
\text{{\scriptsize 0}} & \text{{\scriptsize 0}} & \text{{\scriptsize 0}} & 
\text{{\scriptsize 0}} & \text{{\scriptsize 0}} & \text{{\scriptsize 0}} & 
\text{{\scriptsize 0}} & \text{{\scriptsize 0}} \\ 
\text{{\scriptsize 0}} & \text{{\scriptsize 0}} & \text{{\scriptsize 0}} & 
\text{{\scriptsize 0}} & \text{{\scriptsize 0}} & \text{{\scriptsize 0}} & 
\text{{\scriptsize 0}} & \text{{\scriptsize 0}} \\ 
\text{{\scriptsize 0}} & \text{{\scriptsize 0}} & \text{{\scriptsize 0}} & 
\text{{\scriptsize 0}} & \text{{\scriptsize 0}} & \text{{\scriptsize 0}} & 
\text{{\scriptsize 0}} & \text{{\scriptsize 0}} \\ 
\text{{\scriptsize 0}} & \text{{\scriptsize 0}} & \text{{\scriptsize 0}} & 
\text{{\scriptsize 0}} & \text{{\scriptsize 0}} & \text{{\scriptsize 0}} & 
\text{{\scriptsize 0}} & \text{{\scriptsize 0}} \\ 
\text{{\scriptsize 0}} & \text{{\scriptsize 0}} & \text{{\scriptsize 0}} & 
\text{{\scriptsize 0}} & \text{{\scriptsize 0}} & \text{{\scriptsize 0}} & 
\text{{\scriptsize 0}} & \text{{\scriptsize 0}}
\end{array}
\right) , &  & D_1=\left( 
\begin{array}{cccccccc}
\text{{\scriptsize 1}} & \text{{\scriptsize 0}} & \text{{\scriptsize 1}} & 
\text{{\scriptsize 0}} & \text{{\scriptsize 0}} & \text{{\scriptsize 1}} & 
\text{{\scriptsize 0}} & \text{{\scriptsize 0}} \\ 
\text{{\scriptsize 1}} & \text{{\scriptsize 0}} & \text{{\scriptsize 0}} & 
\text{{\scriptsize 0}} & \text{{\scriptsize 0}} & \text{{\scriptsize 0}} & 
\text{{\scriptsize 0}} & \text{{\scriptsize 0}} \\ 
\text{{\scriptsize 0}} & \text{{\scriptsize 1}} & \text{{\scriptsize 0}} & 
\text{{\scriptsize 1}} & \text{{\scriptsize 1}} & \text{{\scriptsize 0}} & 
\text{{\scriptsize 0}} & \text{{\scriptsize 1}} \\ 
\text{{\scriptsize 0}} & \text{{\scriptsize 1}} & \text{{\scriptsize 0}} & 
\text{{\scriptsize 0}} & \text{{\scriptsize 0}} & \text{{\scriptsize 0}} & 
\text{{\scriptsize 0}} & \text{{\scriptsize 0}} \\ 
\text{{\scriptsize 0}} & \text{{\scriptsize 0}} & \text{{\scriptsize 1}} & 
\text{{\scriptsize 0}} & \text{{\scriptsize 0}} & \text{{\scriptsize 0}} & 
\text{{\scriptsize 0}} & \text{{\scriptsize 1}} \\ 
\text{{\scriptsize 0}} & \text{{\scriptsize 0}} & \text{{\scriptsize 0}} & 
\text{{\scriptsize 1}} & \text{{\scriptsize 0}} & \text{{\scriptsize 0}} & 
\text{{\scriptsize 1}} & \text{{\scriptsize 0}} \\ 
\text{{\scriptsize 0}} & \text{{\scriptsize 0}} & \text{{\scriptsize 0}} & 
\text{{\scriptsize 0}} & \text{{\scriptsize 1}} & \text{{\scriptsize 0}} & 
\text{{\scriptsize 0}} & \text{{\scriptsize 0}} \\ 
\text{{\scriptsize 0}} & \text{{\scriptsize 0}} & \text{{\scriptsize 0}} & 
\text{{\scriptsize 0}} & \text{{\scriptsize 0}} & \text{{\scriptsize 1}} & 
\text{{\scriptsize 1}} & \text{{\scriptsize 0}}
\end{array}
\right)
\end{array}
\]
and $w=(1,2,2,4,3,2,3,3)$.

The eigenvalues of $D_0+D_1$ have minimal polynomial $\xi ^5-3\xi ^4-2\xi
^2-8\xi +8$ and 
\[
\frac{\ln (\limfunc{E}(h_4(N)))}{\ln (n)}\rightarrow
0.7362115557393079316549209... 
\]
as $n\rightarrow \infty $.

\subsection{Sextinomials}

Define $g_5(n)$ to be the number of odd coefficients in $(1+x+\ldots
+x^4+x^5)^n$. Define $P(i,j)$ to be the $(i,j)^{\text{th}}$ entry of
Pascal's sextinomial triangle mod $2$: 
\[
\left( 
\begin{array}{cccccccc}
\QTR{textit}{{\scriptsize {P}}}\text{{\scriptsize (0,0)}} & \QTR{textit}{%
{\scriptsize {P}}}\text{{\scriptsize (0,1)}} & \ldots & \QTR{textit}{%
{\scriptsize {P}}}\text{{\scriptsize (0,4)}} &  &  &  &  \\ 
\QTR{textit}{{\scriptsize {P}}}\text{{\scriptsize (1,0)}} & \QTR{textit}{%
{\scriptsize {P}}}\text{{\scriptsize (1,1)}} & \ldots & \QTR{textit}{%
{\scriptsize {P}}}\text{{\scriptsize (1,4)}} & \QTR{textit}{{\scriptsize {P}}%
}\text{{\scriptsize (1,5)}} & \QTR{textit}{{\scriptsize {P}}}\text{%
{\scriptsize (1,6)}} & {\ldots } & \QTR{textit}{{\scriptsize {P}}}\text{%
{\scriptsize (1,9)}}
\end{array}
\right) =\left( 
\begin{array}{cccccccccc}
\text{{\scriptsize 1}} & \text{{\scriptsize 0}} & \text{{\scriptsize 0}} & 
\text{{\scriptsize 0}} & \text{{\scriptsize 0}} &  &  &  &  &  \\ 
\text{{\scriptsize 1}} & \text{{\scriptsize 1}} & \text{{\scriptsize 1}} & 
\text{{\scriptsize 1}} & \text{{\scriptsize 1}} & \text{{\scriptsize 1}} & 
\text{{\scriptsize 0}} & \text{{\scriptsize 0}} & \text{{\scriptsize 0}} & 
\text{{\scriptsize 0}}
\end{array}
\right) . 
\]
Let $\ell =5$ and $A_0=(1,0,0,0,0)$. We simply summarize: 
\[
A_1=b_0^{1,0}=(1,1,1,0,0); 
\]
\[
\begin{array}{ccccccc}
B_0^{0,0}=A_0, &  & B_0^{1,0}=A_1, &  & B_0^{0,1}=0, &  & B_0^{1,1}=A_1;
\end{array}
\]
\[
A_2=b_1^{0,0}=(1,1,0,0,0); 
\]
\[
\begin{array}{ccccccc}
B_1^{0,0}=A_2, &  & B_1^{1,0}=A_1, &  & B_1^{0,1}=A_0, &  & B_1^{1,1}=A_1;
\end{array}
\]
\[
A_3=b_2^{1,0}=(1,0,0,1,0); 
\]
\[
\begin{array}{ccccccc}
B_2^{0,0}=A_0, &  & B_2^{1,0}=A_3, &  & B_2^{0,1}=A_0, &  & B_2^{1,1}=0;
\end{array}
\]
\[
A_4=b_3^{1,0}=(1,1,0,1,1); 
\]
\[
\begin{array}{ccccccc}
B_3^{0,0}=A_0, &  & B_3^{1,0}=A_4, &  & B_3^{0,1}=A_0, &  & B_3^{1,1}=A_3;
\end{array}
\]
\[
A_5=b_4^{0,0}=(1,0,1,0,0); 
\]
\[
\begin{array}{ccccccc}
B_4^{0,0}=A_5, &  & B_4^{1,0}=A_3, &  & B_4^{0,1}=A_2, &  & B_4^{1,1}=A_3;
\end{array}
\]
\[
\begin{array}{ccccccc}
B_5^{0,0}=A_2, &  & B_5^{1,0}=A_3, &  & B_5^{0,1}=0, &  & B_5^{1,1}=A_3;
\end{array}
\]
\[
\begin{array}{ccc}
D_0=\left( 
\begin{array}{cccccc}
1 & 1 & 2 & 2 & 0 & 0 \\ 
0 & 0 & 0 & 0 & 0 & 0 \\ 
0 & 1 & 0 & 0 & 1 & 1 \\ 
0 & 0 & 0 & 0 & 0 & 0 \\ 
0 & 0 & 0 & 0 & 0 & 0 \\ 
0 & 0 & 0 & 0 & 1 & 0
\end{array}
\right) , &  & D_1=\left( 
\begin{array}{cccccc}
0 & 0 & 0 & 0 & 0 & 0 \\ 
2 & 2 & 0 & 0 & 0 & 0 \\ 
0 & 0 & 0 & 0 & 0 & 0 \\ 
0 & 0 & 1 & 1 & 2 & 2 \\ 
0 & 0 & 0 & 1 & 0 & 0 \\ 
0 & 0 & 0 & 0 & 0 & 0
\end{array}
\right)
\end{array}
\]
and $w=(1,3,2,2,4,2)$.

The eigenvalues of $D_0+D_1$ have minimal polynomial $\xi ^6-4\xi ^5+\xi
^4-\xi ^3+8\xi ^2+11\xi +8$ and 
\[
\frac{\ln (\limfunc{E}(g_5(N)))}{\ln (n)}\rightarrow
0.8194694621655401465959376... 
\]
as $n\rightarrow \infty $.

\subsection{Septinomials}

Define $g_6(n)$ to be the number of odd coefficients in $(1+x+\cdots
+x^5+x^6)^n$. Define $P(i,j)$ to be the $(i,j)^{\text{th}}$ entry of
Pascal's septinomial triangle mod $2$: 
\[
\left( 
\begin{array}{cccccccc}
\QTR{textit}{{\scriptsize {P}}}\text{{\scriptsize (0,0)}} & \QTR{textit}{%
{\scriptsize {P}}}\text{{\scriptsize (0,1)}} & \ldots & \QTR{textit}{%
{\scriptsize {P}}}\text{{\scriptsize (0,5)}} &  &  &  &  \\ 
\QTR{textit}{{\scriptsize {P}}}\text{{\scriptsize (1,0)}} & \QTR{textit}{%
{\scriptsize {P}}}\text{{\scriptsize (1,1)}} & \ldots & \QTR{textit}{%
{\scriptsize {P}}}\text{{\scriptsize (1,5)}} & \QTR{textit}{{\scriptsize {P}}%
}\text{{\scriptsize (1,6)}} & \QTR{textit}{{\scriptsize {P}}}\text{%
{\scriptsize (1,7)}} & {\ldots } & \QTR{textit}{{\scriptsize {P}}}\text{%
{\scriptsize (1,11)}}
\end{array}
\right) =\left( 
\begin{array}{cccccccccccc}
\text{{\scriptsize 1}} & \text{{\scriptsize 0}} & \text{{\scriptsize 0}} & 
\text{{\scriptsize 0}} & \text{{\scriptsize 0}} & \text{{\scriptsize 0}} & 
&  &  &  &  &  \\ 
\text{{\scriptsize 1}} & \text{{\scriptsize 1}} & \text{{\scriptsize 1}} & 
\text{{\scriptsize 1}} & \text{{\scriptsize 1}} & \text{{\scriptsize 1}} & 
\text{{\scriptsize 1}} & \text{{\scriptsize 0}} & \text{{\scriptsize 0}} & 
\text{{\scriptsize 0}} & \text{{\scriptsize 0}} & \text{{\scriptsize 0}}
\end{array}
\right) . 
\]
Let $\ell =6$ and $A_0=(1,0,0,0,0,0)$. We simply summarize: 
\[
\begin{array}{ccc}
A_1=b_0^{1,0}=(1,1,1,1,0,0), &  & A_2=b_0^{1,1}=(1,1,1,0,0,0);
\end{array}
\]
\[
\begin{array}{ccccccc}
B_0^{0,0}=A_0, &  & B_0^{1,0}=A_1, &  & B_0^{0,1}=0, &  & B_0^{1,1}=A_2;
\end{array}
\]
\[
A_3=b_1^{0,0}=(1,1,0,0,0,0); 
\]
\[
\begin{array}{ccccccc}
B_1^{0,0}=A_3, &  & B_1^{1,0}=A_3, &  & B_1^{0,1}=A_3, &  & B_1^{1,1}=A_3;
\end{array}
\]
\[
A_4=b_2^{1,0}=(1,1,1,1,1,0); 
\]
\[
\begin{array}{ccccccc}
B_2^{0,0}=A_3, &  & B_2^{1,0}=A_4, &  & B_2^{0,1}=A_0, &  & B_2^{1,1}=A_3;
\end{array}
\]
\[
\begin{array}{ccccccc}
B_3^{0,0}=A_0, &  & B_3^{1,0}=A_0, &  & B_3^{0,1}=A_0, &  & B_3^{1,1}=A_0;
\end{array}
\]
\[
A_5=b_4^{1,0}=(1,1,1,1,1,1); 
\]
\[
\begin{array}{ccccccc}
B_4^{0,0}=A_2, &  & B_4^{1,0}=A_5, &  & B_4^{0,1}=A_3, &  & B_4^{1,1}=A_0;
\end{array}
\]
\[
\begin{array}{ccccccc}
B_5^{0,0}=A_2, &  & B_5^{1,0}=A_2, &  & B_5^{0,1}=A_2, &  & B_5^{1,1}=A_2;
\end{array}
\]
\[
\begin{array}{ccc}
D_0=\left( 
\begin{array}{cccccc}
1 & 0 & 1 & 2 & 0 & 0 \\ 
0 & 0 & 0 & 0 & 0 & 0 \\ 
0 & 0 & 0 & 0 & 1 & 2 \\ 
0 & 2 & 1 & 0 & 1 & 0 \\ 
0 & 0 & 0 & 0 & 0 & 0 \\ 
0 & 0 & 0 & 0 & 0 & 0
\end{array}
\right) , &  & D_1=\left( 
\begin{array}{cccccc}
0 & 0 & 0 & 2 & 1 & 0 \\ 
1 & 0 & 0 & 0 & 0 & 0 \\ 
1 & 0 & 0 & 0 & 0 & 2 \\ 
0 & 2 & 1 & 0 & 0 & 0 \\ 
0 & 0 & 1 & 0 & 0 & 0 \\ 
0 & 0 & 0 & 0 & 1 & 0
\end{array}
\right)
\end{array}
\]
and $w=(1,4,3,2,5,6)$.

The eigenvalues of $D_0+D_1$ have minimal polynomial $\xi ^6-\xi ^5-2\xi
^4-28\xi ^3+16\xi +64$ and 
\[
\frac{\ln (\limfunc{E}(g_6(N)))}{\ln (n)}\rightarrow
0.8317963967344406899938931... 
\]
as $n\rightarrow \infty $.

\subsection{\label{Lyap}Lyapunov Exponents}

All $D_0$ matrices shown in this paper satisfy $\limfunc{rank}(D_0^q)=1$ for
some positive integer $q$; further, there is a change of coordinates under
which $D_0^q$ is transformed into the matrix whose $(0,0)^{\text{th}}$ entry
is $1$ and all of whose other entries are $0$. That is, there is an
invertible $m\times m$ matrix $Q$ with $Q^{-1}D_0^qQ=e_0e_0^T$. Define 
\[
\begin{array}{ccc}
D_0^{\prime }=Q^{-1}D_0Q, &  & D_1^{\prime }=Q^{-1}D_1Q
\end{array}
\]
and let $z_0$, $z_1$, $\ldots $, $z_{k-2}$, $z_{k-1}$ denote a sequence of
independent random coin tosses (heads=1 and tails=0 with equal probability).
The Lyapunov exponent corresponding to random products of $D_0^{\prime }$
and $D_1^{\prime }$: 
\[
\lambda =\lim_{k\rightarrow \infty }\frac 1k\ln \left\| D_{z_0}^{\prime
}D_{z_1}^{\prime }\cdots D_{z_{k-2}}^{\prime }D_{z_{k-1}}^{\prime }\right\| 
\]
exists almost surely. It turns out that $\lambda /\ln (2)$ is precisely what
we seek to characterize typical growth rates of $g_i(N)$ and $h_j(N)$.

Let $\chi (0^q)$ denote the set of all finite binary words $z$ with no
subwords $0^q$ and with rightmost digit $1$. ($0^1$ means $0$; $0^2$ means $%
00$; $0^3$ means $000$.) Let $\ell (z)$ denote the length of $z$. For
example, $1111$ is the only word of length 4 in $\chi (0)$; $0101$, $0111$, $%
1011$, $1101$, $1111$ are the only words of length 4 in $\chi (00)$; the set 
$\chi (000)$ additionally contains $0011$ and $1001$. It is natural to sort
the elements of $\chi (0^q)$ in terms of increasing length. Write $%
z=z_0z_1\ldots z_{k-2}z_{k-1}$ and $D_z^{\prime }=D_{z_0}^{\prime
}D_{z_1}^{\prime }\cdots D_{z_{k-2}}^{\prime }D_{z_{k-1}}^{\prime }$. By $%
D_z^{\prime }(0,0)$ is meant the upper left corner entry of $D_z^{\prime }$.

Extending earlier work by Pincus \cite{Pi} and Lima \&\ Rahibe \cite{LR},
Moshe \cite{M6} proved that 
\[
\lambda =\frac 1{2^{q+1}(2^q-1)}\dsum\limits_{z\in \chi (0^q)}\frac
1{2^{\ell (z)}}\ln \left| D_z^{\prime }(0,0)\right| . 
\]
This series is attractive, but computationally difficult since the number of
words in $\chi (0^q)$ of length $k$ grows exponentially with increasing $k$.
Summation of the series, coupled with Wynn's $\varepsilon $-process for
accelerating convergence, serves as our primary method for calculating $%
\lambda $. Our secondary method is based on the cycle expansion method
applied to a corresponding Ruelle dynamical zeta function \cite{M1, M2, M3,
Ba}.

In the case of the first trinomial $(1+x+x^2)^n$, we have $q=1$, 
\[
\begin{array}{ccccc}
Q=\left( 
\begin{array}{cc}
1 & -2 \\ 
0 & 1
\end{array}
\right) , &  & D_0^{\prime }=\left( 
\begin{array}{cc}
1 & 0 \\ 
0 & 0
\end{array}
\right) , &  & D_1^{\prime }=\left( 
\begin{array}{cc}
3 & -4 \\ 
1 & -2
\end{array}
\right)
\end{array}
\]
hence \cite{M5} 
\begin{eqnarray*}
\lambda &=&\frac 14\dsum\limits_{k=1}^\infty \frac 1{2^k}\ln \left|
(D_1^{\prime })^k(0,0)\right| =\frac 14\dsum\limits_{k=1}^\infty \frac
1{2^k}\ln \left( \frac{2^{k+2}-(-1)^k}3\right) \\
\ &=&0.4299474333424527201146970... \\
&=&\ln (1.5371767171823579495901403...)
\end{eqnarray*}
and 
\[
\frac{\limfunc{E}(\ln (g_2(N)))}{\ln (n)}\rightarrow \frac \lambda {\ln
(2)}=0.6202830299260946960737425... 
\]
as $n\rightarrow \infty $. This is smaller than $\ln (\varphi )/\ln
(2)=0.694...$, as discussed earlier.

In the case of the quadrinomial $(1+x+x^2+x^3)^n$, we have $q=2$, 
\[
\begin{array}{ccccc}
Q=\left( 
\begin{array}{ccc}
1 & -2 & -2 \\ 
0 & 0 & 1 \\ 
0 & 1 & 0
\end{array}
\right) , &  & D_0^{\prime }=\left( 
\begin{array}{ccc}
1 & 0 & 0 \\ 
0 & 0 & 0 \\ 
0 & 1 & 0
\end{array}
\right) , &  & D_1^{\prime }=\left( 
\begin{array}{ccc}
4 & -4 & -6 \\ 
0 & 2 & 1 \\ 
2 & -4 & -4
\end{array}
\right)
\end{array}
\]
hence 
\[
\lambda =\frac 1{24}\dsum\limits_{z\in \chi (00)}\frac 1{2^{\ell (z)}}\ln
\left| D_z^{\prime }(0,0)\right| =0.34657359... 
\]
and $\lambda /\ln (2)=0.49999999...<0.584....$ We conjecture that $\lambda
/\ln (2)$ equals $1/2$ (the typical growth rate for binomials) and prove
this to be true in section [\ref{Prf}]. The cycle expansion converges slowly
in this case, therefore Moshe's technique is helpful here.

In the case of the second trinomial $(1+x+x^3)^n$, we have $q=2$, 
\[
\begin{array}{ccccc}
Q=\left( 
\begin{array}{cccc}
1 & -3 & -2 & 1 \\ 
0 & 0 & 1 & 0 \\ 
0 & 1 & 0 & -1 \\ 
0 & 0 & 0 & 1
\end{array}
\right) , &  & D_0^{\prime }=\left( 
\begin{array}{cccc}
1 & 0 & 0 & 0 \\ 
0 & 0 & 0 & 0 \\ 
0 & 1 & 0 & 0 \\ 
0 & 0 & 0 & 0
\end{array}
\right) , &  & D_1^{\prime }=\left( 
\begin{array}{cccc}
3 & -6 & -2 & 4 \\ 
0 & 1 & 1 & 0 \\ 
1 & -3 & -2 & 2 \\ 
0 & 1 & 0 & 0
\end{array}
\right)
\end{array}
\]
hence $\lambda =0.45454538229305...$ and $\lambda /\ln
(2)=0.65577036889316...<0.727....$ A related example appears in \cite{M6}.

In the case of the quintinomial $(1+x+x^2+x^3+x^4)^n$, we have $q=2$, 
\[
\begin{array}{ccccc}
Q=\left( 
\begin{array}{cccc}
1 & -3 & -2 & 2 \\ 
0 & 1 & 0 & -2 \\ 
0 & 0 & 1 & 0 \\ 
0 & 0 & 0 & 1
\end{array}
\right) , &  & D_0^{\prime }=\left( 
\begin{array}{cccc}
1 & 0 & 0 & 0 \\ 
0 & 0 & 0 & 0 \\ 
0 & 1 & 0 & 0 \\ 
0 & 0 & 0 & 0
\end{array}
\right) , &  & D_1^{\prime }=\left( 
\begin{array}{cccc}
5 & -10 & -8 & 4 \\ 
1 & -1 & -2 & -2 \\ 
1 & -3 & -2 & 4 \\ 
0 & 1 & 0 & -2
\end{array}
\right)
\end{array}
\]
hence $\lambda =0.504253705692...$ and $\lambda /\ln
(2)=0.727484320552...<0.789....$

In the case of the third trinomial $(1+x+x^4)^n$, we have $q=2$, 
\[
Q=\left( \left( 
\begin{array}{cccccccc}
\text{{\scriptsize 1}} & \text{{\scriptsize -2}} & \text{{\scriptsize -2}} & 
\text{{\scriptsize 0}} & \text{{\scriptsize -1}} & \text{{\scriptsize -2}} & 
\text{{\scriptsize -1}} & \text{{\scriptsize -1}} \\ 
\text{{\scriptsize 0}} & \text{{\scriptsize 1}} & \text{{\scriptsize 0}} & 
\text{{\scriptsize -2}} & \text{{\scriptsize -1}} & \text{{\scriptsize 0}} & 
\text{{\scriptsize -1}} & \text{{\scriptsize -1}} \\ 
\text{{\scriptsize 0}} & \text{{\scriptsize 0}} & \text{{\scriptsize 1}} & 
\text{{\scriptsize 0}} & \text{{\scriptsize 0}} & \text{{\scriptsize 0}} & 
\text{{\scriptsize 0}} & \text{{\scriptsize 0}} \\ 
\text{{\scriptsize 0}} & \text{{\scriptsize 0}} & \text{{\scriptsize 0}} & 
\text{{\scriptsize 1}} & \text{{\scriptsize 0}} & \text{{\scriptsize 0}} & 
\text{{\scriptsize 0}} & \text{{\scriptsize 0}} \\ 
\text{{\scriptsize 0}} & \text{{\scriptsize 0}} & \text{{\scriptsize 0}} & 
\text{{\scriptsize 0}} & \text{{\scriptsize 1}} & \text{{\scriptsize 0}} & 
\text{{\scriptsize 0}} & \text{{\scriptsize 0}} \\ 
\text{{\scriptsize 0}} & \text{{\scriptsize 0}} & \text{{\scriptsize 0}} & 
\text{{\scriptsize 0}} & \text{{\scriptsize 0}} & \text{{\scriptsize 1}} & 
\text{{\scriptsize 0}} & \text{{\scriptsize 0}} \\ 
\text{{\scriptsize 0}} & \text{{\scriptsize 0}} & \text{{\scriptsize 0}} & 
\text{{\scriptsize 0}} & \text{{\scriptsize 0}} & \text{{\scriptsize 0}} & 
\text{{\scriptsize 1}} & \text{{\scriptsize 0}} \\ 
\text{{\scriptsize 0}} & \text{{\scriptsize 0}} & \text{{\scriptsize 0}} & 
\text{{\scriptsize 0}} & \text{{\scriptsize 0}} & \text{{\scriptsize 0}} & 
\text{{\scriptsize 0}} & \text{{\scriptsize 1}}
\end{array}
\right) \right) , 
\]
\[
\begin{array}{ccc}
D_0^{\prime }=\left( 
\begin{array}{cccccccc}
\text{{\scriptsize 1}} & \text{{\scriptsize 0}} & \text{{\scriptsize 0}} & 
\text{{\scriptsize 0}} & \text{{\scriptsize 0}} & \text{{\scriptsize 0}} & 
\text{{\scriptsize 0}} & \text{{\scriptsize 0}} \\ 
\text{{\scriptsize 0}} & \text{{\scriptsize 0}} & \text{{\scriptsize 0}} & 
\text{{\scriptsize 0}} & \text{{\scriptsize 0}} & \text{{\scriptsize 0}} & 
\text{{\scriptsize 0}} & \text{{\scriptsize 0}} \\ 
\text{{\scriptsize 0}} & \text{{\scriptsize 1}} & \text{{\scriptsize 0}} & 
\text{{\scriptsize 0}} & \text{{\scriptsize 0}} & \text{{\scriptsize 0}} & 
\text{{\scriptsize 0}} & \text{{\scriptsize 0}} \\ 
\text{{\scriptsize 0}} & \text{{\scriptsize 0}} & \text{{\scriptsize 0}} & 
\text{{\scriptsize 0}} & \text{{\scriptsize 0}} & \text{{\scriptsize 0}} & 
\text{{\scriptsize 0}} & \text{{\scriptsize 0}} \\ 
\text{{\scriptsize 0}} & \text{{\scriptsize 0}} & \text{{\scriptsize 0}} & 
\text{{\scriptsize 0}} & \text{{\scriptsize 0}} & \text{{\scriptsize 0}} & 
\text{{\scriptsize 0}} & \text{{\scriptsize 0}} \\ 
\text{{\scriptsize 0}} & \text{{\scriptsize 0}} & \text{{\scriptsize 0}} & 
\text{{\scriptsize 0}} & \text{{\scriptsize 0}} & \text{{\scriptsize 0}} & 
\text{{\scriptsize 0}} & \text{{\scriptsize 0}} \\ 
\text{{\scriptsize 0}} & \text{{\scriptsize 0}} & \text{{\scriptsize 0}} & 
\text{{\scriptsize 0}} & \text{{\scriptsize 0}} & \text{{\scriptsize 0}} & 
\text{{\scriptsize 0}} & \text{{\scriptsize 0}} \\ 
\text{{\scriptsize 0}} & \text{{\scriptsize 0}} & \text{{\scriptsize 0}} & 
\text{{\scriptsize 0}} & \text{{\scriptsize 0}} & \text{{\scriptsize 0}} & 
\text{{\scriptsize 0}} & \text{{\scriptsize 0}}
\end{array}
\right) , &  & D_1^{\prime }=\left( 
\begin{array}{cccccccc}
\text{{\scriptsize 3}} & \text{{\scriptsize 0}} & \text{{\scriptsize -2}} & 
\text{{\scriptsize -8}} & \text{{\scriptsize -4}} & \text{{\scriptsize -2}}
& \text{{\scriptsize -4}} & \text{{\scriptsize -4}} \\ 
\text{{\scriptsize 1}} & \text{{\scriptsize 0}} & \text{{\scriptsize -1}} & 
\text{{\scriptsize -4}} & \text{{\scriptsize -2}} & \text{{\scriptsize -1}}
& \text{{\scriptsize -2}} & \text{{\scriptsize -2}} \\ 
\text{{\scriptsize 0}} & \text{{\scriptsize 1}} & \text{{\scriptsize 0}} & 
\text{{\scriptsize -1}} & \text{{\scriptsize 0}} & \text{{\scriptsize 0}} & 
\text{{\scriptsize -1}} & \text{{\scriptsize 0}} \\ 
\text{{\scriptsize 0}} & \text{{\scriptsize 1}} & \text{{\scriptsize 0}} & 
\text{{\scriptsize -2}} & \text{{\scriptsize -1}} & \text{{\scriptsize 0}} & 
\text{{\scriptsize -1}} & \text{{\scriptsize -1}} \\ 
\text{{\scriptsize 0}} & \text{{\scriptsize 0}} & \text{{\scriptsize 1}} & 
\text{{\scriptsize 0}} & \text{{\scriptsize 0}} & \text{{\scriptsize 0}} & 
\text{{\scriptsize 0}} & \text{{\scriptsize 1}} \\ 
\text{{\scriptsize 0}} & \text{{\scriptsize 0}} & \text{{\scriptsize 0}} & 
\text{{\scriptsize 1}} & \text{{\scriptsize 0}} & \text{{\scriptsize 0}} & 
\text{{\scriptsize 1}} & \text{{\scriptsize 0}} \\ 
\text{{\scriptsize 0}} & \text{{\scriptsize 0}} & \text{{\scriptsize 0}} & 
\text{{\scriptsize 0}} & \text{{\scriptsize 1}} & \text{{\scriptsize 0}} & 
\text{{\scriptsize 0}} & \text{{\scriptsize 0}} \\ 
\text{{\scriptsize 0}} & \text{{\scriptsize 0}} & \text{{\scriptsize 0}} & 
\text{{\scriptsize 0}} & \text{{\scriptsize 0}} & \text{{\scriptsize 1}} & 
\text{{\scriptsize 1}} & \text{{\scriptsize 0}}
\end{array}
\right)
\end{array}
\]
hence $\lambda =0.45759385431410...$ and $\lambda /\ln
(2)=0.66016838436022...<0.736....$

In the case of the sextinomial $(1+x+\cdots +x^4+x^5)^n$, we have $q=3$, 
\[
\begin{array}{ccccc}
Q=\left( 
\begin{array}{cccccc}
\text{{\scriptsize 1}} & \text{{\scriptsize -1}} & \text{{\scriptsize -2}} & 
\text{{\scriptsize -2}} & \text{{\scriptsize -1}} & \text{{\scriptsize -2}}
\\ 
\text{{\scriptsize 0}} & \text{{\scriptsize -1}} & \text{{\scriptsize 0}} & 
\text{{\scriptsize 0}} & \text{{\scriptsize 1}} & \text{{\scriptsize 0}} \\ 
\text{{\scriptsize 0}} & \text{{\scriptsize 0}} & \text{{\scriptsize 0}} & 
\text{{\scriptsize 1}} & \text{{\scriptsize 0}} & \text{{\scriptsize 0}} \\ 
\text{{\scriptsize 0}} & \text{{\scriptsize 0}} & \text{{\scriptsize 0}} & 
\text{{\scriptsize 0}} & \text{{\scriptsize 0}} & \text{{\scriptsize 1}} \\ 
\text{{\scriptsize 0}} & \text{{\scriptsize 1}} & \text{{\scriptsize 0}} & 
\text{{\scriptsize 0}} & \text{{\scriptsize 0}} & \text{{\scriptsize 0}} \\ 
\text{{\scriptsize 0}} & \text{{\scriptsize 0}} & \text{{\scriptsize 1}} & 
\text{{\scriptsize 0}} & \text{{\scriptsize -1}} & \text{{\scriptsize 0}}
\end{array}
\right) , &  & D_0^{\prime }=\left( 
\begin{array}{cccccc}
\text{{\scriptsize 1}} & \text{{\scriptsize 0}} & \text{{\scriptsize 0}} & 
\text{{\scriptsize 0}} & \text{{\scriptsize 0}} & \text{{\scriptsize 0}} \\ 
\text{{\scriptsize 0}} & \text{{\scriptsize 0}} & \text{{\scriptsize 0}} & 
\text{{\scriptsize 0}} & \text{{\scriptsize 0}} & \text{{\scriptsize 0}} \\ 
\text{{\scriptsize 0}} & \text{{\scriptsize 1}} & \text{{\scriptsize 0}} & 
\text{{\scriptsize 0}} & \text{{\scriptsize 0}} & \text{{\scriptsize 0}} \\ 
\text{{\scriptsize 0}} & \text{{\scriptsize 0}} & \text{{\scriptsize 1}} & 
\text{{\scriptsize 0}} & \text{{\scriptsize 0}} & \text{{\scriptsize 0}} \\ 
\text{{\scriptsize 0}} & \text{{\scriptsize 0}} & \text{{\scriptsize 0}} & 
\text{{\scriptsize 0}} & \text{{\scriptsize 0}} & \text{{\scriptsize 0}} \\ 
\text{{\scriptsize 0}} & \text{{\scriptsize 0}} & \text{{\scriptsize 0}} & 
\text{{\scriptsize 0}} & \text{{\scriptsize 0}} & \text{{\scriptsize 0}}
\end{array}
\right) , &  & D_1^{\prime }=\left( 
\begin{array}{cccccc}
\text{{\scriptsize 6}} & \text{{\scriptsize -8}} & \text{{\scriptsize -8}} & 
\text{{\scriptsize -10}} & \text{{\scriptsize -4}} & \text{{\scriptsize -6}}
\\ 
\text{{\scriptsize 0}} & \text{{\scriptsize 0}} & \text{{\scriptsize 0}} & 
\text{{\scriptsize 0}} & \text{{\scriptsize 0}} & \text{{\scriptsize 1}} \\ 
\text{{\scriptsize 2}} & \text{{\scriptsize -4}} & \text{{\scriptsize -4}} & 
\text{{\scriptsize -4}} & \text{{\scriptsize 0}} & \text{{\scriptsize -3}}
\\ 
\text{{\scriptsize 0}} & \text{{\scriptsize 0}} & \text{{\scriptsize 0}} & 
\text{{\scriptsize 0}} & \text{{\scriptsize 0}} & \text{{\scriptsize 0}} \\ 
\text{{\scriptsize 2}} & \text{{\scriptsize -4}} & \text{{\scriptsize -4}} & 
\text{{\scriptsize -4}} & \text{{\scriptsize 0}} & \text{{\scriptsize -3}}
\\ 
\text{{\scriptsize 0}} & \text{{\scriptsize 2}} & \text{{\scriptsize 2}} & 
\text{{\scriptsize 1}} & \text{{\scriptsize -2}} & \text{{\scriptsize 1}}
\end{array}
\right)
\end{array}
\]
hence 
\[
\lambda =\frac 1{112}\dsum\limits_{z\in \chi (000)}\frac 1{2^{\ell (z)}}\ln
\left| D_z^{\prime }(0,0)\right| =0.5344481528... 
\]
and $\lambda /\ln (2)=0.7710456996...<0.819....$

In the case of the septinomial $(1+x+\cdots +x^5+x^6)^n$, we have $q=3$, 
\[
\begin{array}{ccccc}
Q=\left( 
\begin{array}{cccccc}
\text{{\scriptsize 1}} & \text{{\scriptsize -6}} & \text{{\scriptsize -6}} & 
\text{{\scriptsize -4}} & \text{{\scriptsize 1}} & \text{{\scriptsize 0}} \\ 
\text{{\scriptsize 0}} & \text{{\scriptsize 0}} & \text{{\scriptsize 0}} & 
\text{{\scriptsize 0}} & \text{{\scriptsize 0}} & \text{{\scriptsize 1}} \\ 
\text{{\scriptsize 0}} & \text{{\scriptsize 0}} & \text{{\scriptsize 2}} & 
\text{{\scriptsize 0}} & \text{{\scriptsize -1}} & \text{{\scriptsize 0}} \\ 
\text{{\scriptsize 0}} & \text{{\scriptsize 0}} & \text{{\scriptsize 0}} & 
\text{{\scriptsize 2}} & \text{{\scriptsize 0}} & \text{{\scriptsize 0}} \\ 
\text{{\scriptsize 0}} & \text{{\scriptsize 0}} & \text{{\scriptsize 0}} & 
\text{{\scriptsize 0}} & \text{{\scriptsize 1}} & \text{{\scriptsize -2}} \\ 
\text{{\scriptsize 0}} & \text{{\scriptsize 1}} & \text{{\scriptsize 0}} & 
\text{{\scriptsize 0}} & \text{{\scriptsize -}}\tfrac 12 & \text{%
{\scriptsize 1}}
\end{array}
\right) , &  & D_0^{\prime }=\left( 
\begin{array}{cccccc}
\text{{\scriptsize 1}} & \text{{\scriptsize 0}} & \text{{\scriptsize 0}} & 
\text{{\scriptsize 0}} & \text{{\scriptsize 0}} & \text{{\scriptsize 0}} \\ 
\text{{\scriptsize 0}} & \text{{\scriptsize 0}} & \text{{\scriptsize 0}} & 
\text{{\scriptsize 0}} & \text{{\scriptsize 0}} & \text{{\scriptsize 0}} \\ 
\text{{\scriptsize 0}} & \text{{\scriptsize 1}} & \text{{\scriptsize 0}} & 
\text{{\scriptsize 0}} & \text{{\scriptsize 0}} & \text{{\scriptsize 0}} \\ 
\text{{\scriptsize 0}} & \text{{\scriptsize 0}} & \text{{\scriptsize 1}} & 
\text{{\scriptsize 0}} & \text{{\scriptsize 0}} & \text{{\scriptsize 0}} \\ 
\text{{\scriptsize 0}} & \text{{\scriptsize 0}} & \text{{\scriptsize 0}} & 
\text{{\scriptsize 0}} & \text{{\scriptsize 0}} & \text{{\scriptsize 0}} \\ 
\text{{\scriptsize 0}} & \text{{\scriptsize 0}} & \text{{\scriptsize 0}} & 
\text{{\scriptsize 0}} & \text{{\scriptsize 0}} & \text{{\scriptsize 0}}
\end{array}
\right) , &  & D_1^{\prime }=\left( 
\begin{array}{cccccc}
\text{{\scriptsize 7}} & \text{{\scriptsize -36}} & \text{{\scriptsize -28}}
& \text{{\scriptsize -24}} & \text{{\scriptsize 4}} & \text{{\scriptsize -4}}
\\ 
\text{{\scriptsize 0}} & \text{{\scriptsize 0}} & \text{{\scriptsize 1}} & 
\text{{\scriptsize 0}} & \tfrac 12 & \text{{\scriptsize -2}} \\ 
\tfrac 32 & \text{{\scriptsize -8}} & \text{{\scriptsize -8}} & \text{%
{\scriptsize -6}} & \tfrac 12 & \text{{\scriptsize 1}} \\ 
\text{{\scriptsize 0}} & \text{{\scriptsize 0}} & \text{{\scriptsize 1}} & 
\text{{\scriptsize 0}} & \text{{\scriptsize -}}\tfrac 12 & \text{%
{\scriptsize 1}} \\ 
\text{{\scriptsize 2}} & \text{{\scriptsize -12}} & \text{{\scriptsize -10}}
& \text{{\scriptsize -8}} & \text{{\scriptsize 1}} & \text{{\scriptsize 0}}
\\ 
\text{{\scriptsize 1}} & \text{{\scriptsize -6}} & \text{{\scriptsize -6}} & 
\text{{\scriptsize -4}} & \text{{\scriptsize 1}} & \text{{\scriptsize 0}}
\end{array}
\right)
\end{array}
\]
hence $\lambda =0.53765282...$ and $\lambda /\ln (2)=0.77566905...<0.831....$

Let us conclude by mentioning relevant sequences in Sloane's online
encyclopedia. A001316 is $f(n)$; A000120 is $\ln (f(n))/\ln (2)$; A006046 is 
$\sum_{k<n}f(k)$; A000788 is $\sum_{k<n}\ln (f(k))/\ln (2)$. A071053 is $%
g(n) $; A134659 is $\sum_{k<n}g(k)$. A134660 is $g_3(n)$; A036555 is $\ln
(g_3(n))/\ln (2)$. A134661 is $h_3(n)$. A134662 is $h_4(n)$. A007318,
A027907, A008287, A035343, A063260, A063265 are Pascal's triangles
associated with $(1+x+\cdots +x^{r-1}+x^r)^n$ for $r=1,\ldots ,6$; A038717
and\ A134663 are likewise for $(1+x+x^3)^n$ and\ $(1+x+x^4)^n$. It is
well-known that $\ln (f(n))/\ln (2)$ is the number of $1$s in the binary
expansion of $n$, but scarcely noticed that $\ln (g_3(n))/\ln (2)$ is the
number of $1$s in the binary expansion of $3n$.

\subsection{Pascal's Rhombus}

The sequence of polynomials giving Pascal's trinomial triangle arises from
the first-order recurrence 
\[
\begin{array}{ccc}
p_n(x)=(1+x+x^2)p_{n-1}(x), &  & p_0(x)=1.
\end{array}
\]
Pascal's rhombus \cite{GKMT}, by contrast, arises from the second-order
recurrence 
\[
\begin{array}{ccccc}
p_n(x)=(1+x+x^2)p_{n-1}(x)+x^2p_{n-2}(x), &  & p_1(x)=1+x+x^2, &  & p_0(x)=1.
\end{array}
\]
In this addendum, we perform the same analysis as in the preceding. Define $%
u(n)$ to be the number of odd coefficients in $p_n(x)$. Define $P(i,j)$ to
be the $(i,j)^{\text{th}}$ entry of the associated ``rhombus'' mod $2$: 
\[
\left( 
\begin{array}{cccccccc}
\QTR{textit}{{\scriptsize {P}}}\text{{\scriptsize (0,0)}} & \QTR{textit}{%
{\scriptsize {P}}}\text{{\scriptsize (0,1)}} & \QTR{textit}{{\scriptsize {P}}%
}\text{{\scriptsize (0,2)}} & \QTR{textit}{{\scriptsize {P}}}\text{%
{\scriptsize (0,3)}} &  &  &  &  \\ 
\QTR{textit}{{\scriptsize {P}}}\text{{\scriptsize (1,0)}} & \QTR{textit}{%
{\scriptsize {P}}}\text{{\scriptsize (1,1)}} & \QTR{textit}{{\scriptsize {P}}%
}\text{{\scriptsize (1,2)}} & \QTR{textit}{{\scriptsize {P}}}\text{%
{\scriptsize (1,3)}} & \QTR{textit}{{\scriptsize {P}}}\text{{\scriptsize %
(1,4)}} & \QTR{textit}{{\scriptsize {P}}}\text{{\scriptsize (1,5)}} &  &  \\ 
\QTR{textit}{{\scriptsize {P}}}\text{{\scriptsize (2,0)}} & \QTR{textit}{%
{\scriptsize {P}}}\text{{\scriptsize (2,1)}} & \QTR{textit}{{\scriptsize {P}}%
}\text{{\scriptsize (2,2)}} & \QTR{textit}{{\scriptsize {P}}}\text{%
{\scriptsize (2,3)}} & \QTR{textit}{{\scriptsize {P}}}\text{{\scriptsize %
(2,4)}} & \QTR{textit}{{\scriptsize {P}}}\text{{\scriptsize (2,5)}} & 
\QTR{textit}{{\scriptsize {P}}}\text{{\scriptsize (2,6)}} & \QTR{textit}{%
{\scriptsize {P}}}\text{{\scriptsize (2,7)}}
\end{array}
\right) =\left( 
\begin{array}{cccccccc}
\text{{\scriptsize 1}} & \text{{\scriptsize 0}} & \text{{\scriptsize 0}} & 
\text{{\scriptsize 0}} &  &  &  &  \\ 
\text{{\scriptsize 1}} & \text{{\scriptsize 1}} & \text{{\scriptsize 1}} & 
\text{{\scriptsize 0}} & \text{{\scriptsize 0}} & \text{{\scriptsize 0}} & 
&  \\ 
\text{{\scriptsize 1}} & \text{{\scriptsize 0}} & \text{{\scriptsize 0}} & 
\text{{\scriptsize 0}} & \text{{\scriptsize 1}} & \text{{\scriptsize 0}} & 
\text{{\scriptsize 0}} & \text{{\scriptsize 0}}
\end{array}
\right) . 
\]
The objects $A_i$ and $B_i$ are no longer $\ell $-vectors but $2\times \ell $
matrices. It is mandatory that $A_i(0,0)=1$ or $A_i(1,0)=1$ (or both) for
every $i$. Let $\ell =4$ and 
\[
A_0=\left( 
\begin{array}{cccc}
1 & 0 & 0 & 0 \\ 
1 & 1 & 1 & 0
\end{array}
\right) . 
\]
We simply summarize: 
\[
\begin{array}{ccc}
A_1=b_0^{0,0}=\left( 
\begin{array}{cccc}
1 & 0 & 0 & 0 \\ 
1 & 0 & 1 & 0
\end{array}
\right) , &  & A_2=b_0^{1,0}=\left( 
\begin{array}{cccc}
1 & 1 & 0 & 0 \\ 
1 & 0 & 0 & 1
\end{array}
\right) ;
\end{array}
\]
\[
\begin{array}{ccccccc}
B_0^{0,0}=A_1, &  & B_0^{1,0}=A_2, &  & B_0^{0,1}=0, &  & B_0^{1,1}=A_0;
\end{array}
\]
\[
\begin{array}{ccc}
A_3=b_1^{0,1}=\left( 
\begin{array}{cccc}
0 & 0 & 0 & 0 \\ 
1 & 1 & 0 & 0
\end{array}
\right) , &  & A_4=B_1^{1,1}=\left( 
\begin{array}{cccc}
0 & 0 & 0 & 0 \\ 
1 & 0 & 0 & 0
\end{array}
\right) ;
\end{array}
\]
\[
\begin{array}{ccccccc}
B_1^{0,0}=A_0, &  & B_1^{1,0}=A_2, &  & B_1^{0,1}=A_3, &  & B_1^{1,1}=A_4;
\end{array}
\]
\[
\begin{array}{ccccccc}
B_2^{0,0}=A_1, &  & B_2^{1,0}=A_0, &  & B_2^{0,1}=A_1, &  & B_2^{1,1}=A_0;
\end{array}
\]
\[
\begin{array}{ccccccc}
B_3^{0,0}=A_4, &  & B_3^{1,0}=A_1, &  & B_3^{0,1}=A_4, &  & B_3^{1,1}=A_1;
\end{array}
\]
\[
\begin{array}{ccccccc}
B_4^{0,0}=A_3, &  & B_4^{1,0}=A_1, &  & B_4^{0,1}=A_4, &  & B_4^{1,1}=0;
\end{array}
\]
\[
\begin{array}{ccc}
D_0=\left( 
\begin{array}{ccccc}
0 & 1 & 0 & 0 & 0 \\ 
1 & 0 & 2 & 0 & 0 \\ 
0 & 0 & 0 & 0 & 0 \\ 
0 & 1 & 0 & 0 & 1 \\ 
0 & 0 & 0 & 2 & 1
\end{array}
\right) , &  & D_1=\left( 
\begin{array}{ccccc}
1 & 0 & 2 & 0 & 0 \\ 
0 & 0 & 0 & 2 & 1 \\ 
1 & 1 & 0 & 0 & 0 \\ 
0 & 0 & 0 & 0 & 0 \\ 
0 & 1 & 0 & 0 & 0
\end{array}
\right)
\end{array}
\]
and $w=(1,1,2,0,0)$. The maximal eigenvalue $\mu _1$ of $D_0+D_1$ is $(3+%
\sqrt{17})/2$; therefore 
\[
\frac{\ln (\limfunc{E}(u(N)))}{\ln (n)}\rightarrow \frac{\ln ((3+\sqrt{17}%
)/4)}{\ln (2)}=0.8325063835804514437981667... 
\]
as $n\rightarrow \infty $. To compute $\lim_{n\rightarrow \infty }\limfunc{E}%
(\ln (u(N)))/\ln (n)$ via Moshe's technique requires a binary word $z$ for
which the product $D_z$ satisfies $\limfunc{rank}(D_z)=1$. No such word $z$
exists, therefore our primary method is inapplicable here. Our secondary
method gives $\lambda =0.57331379313...$, hence $\lambda /\ln
(2)=0.82711696622....$

We mention the Fibonacci polynomials 
\[
\begin{array}{ccccc}
p_n(x)=x\,p_{n-1}(x)+p_{n-2}(x), &  & p_1(x)=x, &  & p_0(x)=1
\end{array}
\]
and that the number $v(n)$ of odd coefficients in $p_n(x)$ is the $n^{\text{%
th}}$ term of Stern's sequence \cite{Ur, CW, R1} 
\[
\begin{array}{ccc}
v(2n+1)=v(n), &  & v(2n)=v(n)+v(n-1).
\end{array}
\]
It follows that 
\[
\begin{array}{ccccc}
D_0=\left( 
\begin{array}{cc}
1 & 0 \\ 
1 & 1
\end{array}
\right) , &  & D_1=\left( 
\begin{array}{cc}
1 & 1 \\ 
0 & 1
\end{array}
\right) , &  & w=(1,0);
\end{array}
\]
\[
\begin{array}{ccc}
\dfrac{\ln (\limfunc{E}(v(N)))}{\ln (n)}\rightarrow \dfrac{\ln (3/2)}{\ln (2)%
}, &  & \dfrac{\limfunc{E}(\ln (v(N)))}{\ln (n)}\rightarrow \dfrac \lambda
{\ln (2)}=0.571613901650254...
\end{array}
\]
as $n\rightarrow \infty $. The estimate $\lambda =0.396212564297744...$ is
new (as far as we know), but the challenge of computing $\lambda $ was posed
long ago \cite{CLM}.

Relevant sequences in Sloane's online encyclopedia are A059319 for $u(n)$,
A059317 for Pascal's rhombus, A002487 for $v(n)$ and A049310 for
``Fibonacci's rhombus''. We recall Glaisher's theorem that $f(n)$ is the
number of odd binomial coefficients of the form $\tbinom nm$, $0\leq m\leq n$%
; this is mirrored by Carlitz's theorem that $v(n)$ is the number of odd
binomial coefficients of the form $\tbinom{n-m}m$, $0\leq 2m\leq n$.

\subsection{Extreme Values}

The function $f(n)$ has ``maximum growth'' $\approx n^1$ in the sense that 
\[
\limfunc{limsup}_{n\rightarrow \infty }\frac 1{\ln (n)}\ln (f(n))=1. 
\]
To evaluate the limit superior, we note that 
\[
\begin{array}{ccc}
f(n)<f(2^k-1)=2^k &  & \text{for all }n<2^k-1
\end{array}
\]
for each $k=0,1,2,\ldots $. No other disjoint subsequence possesses this
property and 
\[
\lim_{k\rightarrow \infty }\frac{\ln (2^k)}{\ln (2^k-1)}=1. 
\]
The functions $g(n)$ and $u(n)$ likewise satisfy 
\[
\limfunc{limsup}_{n\rightarrow \infty }\frac 1{\ln (n)}\ln (g(n))=\limfunc{%
limsup}_{n\rightarrow \infty }\frac 1{\ln (n)}\ln (u(n))=1. 
\]
In the case of $g(n)$, we have 
\[
g(n)<\left\{ 
\begin{array}{lll}
g(2^k-1)=\dfrac{2^{k+2}-(-1)^k}3 &  & \text{for all }n<2^k-1, \\ 
g(3\cdot 2^{k+1}-1)=2^{k+3}+(-1)^k &  & \text{for all }n<3\cdot 2^{k+1}-1,
\\ 
g(11\cdot 2^{k+1}-1)=3(2^{k+3}+(-1)^k) &  & \text{for all }n<11\cdot
2^{k+1}-1
\end{array}
\right. 
\]
for each $k=0,1,2,\ldots $. No other disjoint subsequence possesses this
property and 
\[
\lim_{k\rightarrow \infty }\frac{\ln \left( \tfrac{2^{k+2}-(-1)^k}3\right) }{%
\ln (2^k-1)}=\lim_{k\rightarrow \infty }\frac{\ln (2^{k+3}+(-1)^k)}{\ln
(3\cdot 2^{k+1}-1)}=\lim_{k\rightarrow \infty }\frac{\ln (3(2^{k+3}+(-1)^k))%
}{\ln (11\cdot 2^{k+1}-1)}=1. 
\]
In the case of $u(n)$, we have 
\[
u(n)<\left\{ 
\begin{array}{lll}
u(5)=6 &  & \text{for all }n<5, \\ 
u(37)=45 &  & \text{for all }n<37, \\ 
u(2^k-1)=\dfrac{2^{k+2}-(-1)^k}3 &  & \text{for all }n<2^k-1, \\ 
u(5\cdot 2^{k+1}-1)=\dfrac{5\cdot 2^{k+3}+(-1)^k}3 &  & \text{for all }%
n<5\cdot 2^{k+1}-1, \\ 
u(2\cdot 4^{k+2}-7)=\dfrac{5\cdot 4^{k+2}+12k+1}3 &  & \text{for all }%
n<2\cdot 4^{k+2}-7
\end{array}
\right. 
\]
for each $k=0,1,2,\ldots $. No other disjoint subsequence possesses this
property and 
\[
\lim_{k\rightarrow \infty }\frac{\ln \left( \tfrac{2^{k+2}-(-1)^k}3\right) }{%
\ln (2^k-1)}=\lim_{k\rightarrow \infty }\frac{\ln \left( \tfrac{5\cdot
2^{k+3}+(-1)^k}3\right) }{\ln (5\cdot 2^{k+1}-1)}=\lim_{k\rightarrow \infty }%
\frac{\ln \left( \tfrac{5\cdot 4^{k+2}+12k+1}3\right) }{\ln (2\cdot
4^{k+2}-7)}=1. 
\]
By contrast, the function $v(n)$ satisfies \cite{CW} 
\[
\limfunc{limsup}_{n\rightarrow \infty }\frac 1{\ln (n)}\ln (v(n))=\frac{\ln
(\varphi )}{\ln (2)} 
\]
where $\varphi =(1+\sqrt{5})/2$. Let $a_0=0$, $a_1=1$, $a_j=a_{j-1}+a_{j-2}$
denote the Fibonacci sequence. We have 
\[
v(n)<\left\{ 
\begin{array}{lll}
v\left( \dfrac{2(4^k-1)}3\right) =a_{2k+1} &  & \text{for all }n<\dfrac{%
2(4^k-1)}3, \\ 
v\left( \dfrac{4(4^k-1)}3\right) =a_{2k+2} &  & \text{for all }n<\dfrac{%
4(4^k-1)}3
\end{array}
\right. 
\]
for each $k=0,1,2,\ldots $. No other disjoint subsequence possesses this
property and 
\[
\lim_{k\rightarrow \infty }\frac{\ln (a_{2k+1})}{\ln \left( \tfrac{2(4^k-1)}%
3\right) }=\lim_{k\rightarrow \infty }\frac{\ln (a_{2k+2})}{\ln \left( 
\tfrac{4(4^k-1)}3\right) }=\frac{\ln (\varphi )}{\ln (2)}. 
\]
Our analyses in this addendum are aided by the $D_0$, $D_1$, $w$ matrices
from earlier sections. ``Minimum growth'', defined with limit superior
replaced by limit inferior, is not as interesting for three of the cases
since 
\[
\begin{array}{ccccc}
f(2^k)=2, &  & g(2^k)=3, &  & v(2^k-1)=1
\end{array}
\]
always. The remaining case, $u(n)$, resists all attempts at simplification.

\subsection{Variability}

We shall be very brief here. Let $N$ denote a uniform random integer between 
$0$ and $n-1$. Kirschenhofer \cite{Ki} proved that $f(N)$ has ``typical
dispersion'' $\approx n^{\ln (2)/4}$ in the sense that 
\[
\limfunc{Var}(\ln (f(N)))\sim \frac{\ln (2)}4\ln (n)
\]
as $n\rightarrow \infty ;$ equivalently, 
\[
\lim_{n\rightarrow \infty }\frac 1{\ln (n)}\left[ \frac
1n\dsum\limits_{k=0}^{n-1}\ln (f(k))^2-\left( \frac
1n\dsum\limits_{k=0}^{n-1}\ln (f(k))\right) ^2\right] =\frac{\ln (2)}4.
\]
In fact, this is related to a generalized Lyapunov exponent of order two
corresponding to random products of $D_0=(1)$ and $D_1=(2)$. We defer
further study of such quantities to a later paper.

Define $\psi =(5+\sqrt{17})/4$ for convenience. As another example, $v(N)$
has ``average dispersion'' $\approx n^{\ln (\psi )/\ln (2)}$ in the sense
that 
\[
\ln (\limfunc{Var}(v(N)))\sim \frac{\ln (\psi )}{\ln (2)}\ln (n)
\]
as $n\rightarrow \infty ;$ equivalently, 
\[
\lim_{n\rightarrow \infty }\frac 1{\ln (n)}\ln \left[ \frac
1n\dsum\limits_{k=0}^{n-1}v(k)^2-\left( \frac
1n\dsum\limits_{k=0}^{n-1}v(k)\right) ^2\right] =\frac{\ln (\psi )}{\ln (2)}.
\]
More details on the latter result can be found in \cite{R2, MM, C1, C2, C3,
C4, Al}. In fact, $\psi $ is the leading eigenvalue of $\frac 12(D_0\otimes
D_0+D_1\otimes D_1)$, where $\otimes $ is the direct or Kronecker product of
matrices 
\[
\begin{array}{ccc}
D_0=\left( 
\begin{array}{cc}
1 & 0 \\ 
1 & 1
\end{array}
\right) , &  & D_1=\left( 
\begin{array}{cc}
1 & 1 \\ 
0 & 1
\end{array}
\right) ;
\end{array}
\]
thus 
\[
\begin{array}{ccc}
D_0\otimes D_0=\left( 
\begin{array}{cccc}
1 & 0 & 0 & 0 \\ 
1 & 1 & 0 & 0 \\ 
1 & 0 & 1 & 0 \\ 
1 & 1 & 1 & 1
\end{array}
\right) , &  & D_1\otimes D_1=\left( 
\begin{array}{cccc}
1 & 1 & 1 & 1 \\ 
0 & 1 & 0 & 1 \\ 
0 & 0 & 1 & 1 \\ 
0 & 0 & 0 & 1
\end{array}
\right) .
\end{array}
\]
In the same way, 
\[
\begin{array}{ccc}
\dfrac{\ln (\limfunc{Var}(f(N)))}{\ln (n)}\rightarrow \dfrac{\ln (5/2)}{\ln
(2)}, &  & \dfrac{\ln (\limfunc{Var}(g(N)))}{\ln (n)}\rightarrow \dfrac{\ln
(\xi )}{\ln (2)}
\end{array}
\]
where $\xi =2.813...$ has minimal polynomial $\xi ^3-2\xi ^2-3\xi +2$ and 
\[
\dfrac{\ln (\limfunc{Var}(u(N)))}{\ln (n)}\rightarrow \dfrac{\ln (\eta )}{%
\ln (2)}
\]
where $\eta =3.194...$ has minimal polynomial $4\eta ^7-8\eta ^6-25\eta
^5+22\eta ^4+24\eta ^3+16\eta ^2+\eta -2$.

\subsection{\label{Prf}Proof of Conjecture}

Since the Lyapunov exponent is defined as an almost-sure limit: 
\[
\lambda =\lim_{k\rightarrow \infty }\frac 1k\ln \left\| D_{z_0}D_{z_1}\cdots
D_{z_{k-2}}D_{z_{k-1}}\right\| 
\]
we may assume that $z_0=1$. Every binary word $z$ thus looks like 
\[
1\,0^{j_0}\,1\,0^{j_1}\,1\,0^{j_2}\,\ldots \,1\,0^{j_{n-2}}\,1\,0^{j_{n-1}} 
\]
where each $j_i\geq 0$. We introduce a rewording of $D_z$: 
\[
\tilde D_{j_0}\tilde D_{j_1}\tilde D_{j_2}\cdots \tilde D_{j_{n-2}}\tilde
D_{j_{n-1}}=\left( D_1D_0^{j_0}\right) \left( D_1D_0^{j_1}\right) \left(
D_1D_0^{j_2}\right) \cdots \left( D_1D_0^{j_{n-2}}\right) \left(
D_1D_0^{j_{n-1}}\right) 
\]
and note that the probability associated with $\tilde D_{j_i}$ is $%
1/2^{j_i+1}$. Clearly 
\[
\lambda =\lim_{n\rightarrow \infty }\frac 1{2n}\ln \left\| \tilde
D_{j_0}\tilde D_{j_1}\cdots \tilde D_{j_{n-2}}\tilde D_{j_{n-1}}\right\| 
\]
since the mean length of $\tilde D_{j_i}$ is $1+\tsum\nolimits_{j\geq
0}j/2^{j+1}$. Also 
\[
\begin{array}{ccc}
\tilde D_0=D_1=\left( 
\begin{array}{ccc}
0 & 0 & 0 \\ 
2 & 0 & 0 \\ 
0 & 1 & 2
\end{array}
\right) , &  & \tilde D_1=D_1D_0=\left( 
\begin{array}{ccc}
0 & 0 & 0 \\ 
2 & 4 & 0 \\ 
0 & 0 & 1
\end{array}
\right) ,
\end{array}
\]
\[
\tilde D_j=D_1D_0^j=\left( 
\begin{array}{ccc}
0 & 0 & 0 \\ 
2 & 4 & 4 \\ 
0 & 0 & 0
\end{array}
\right) 
\]
for all $j\geq 2$. Hence the rewording actually consists of only three
matrices $\tilde D_0$, $\tilde D_1$, $\tilde D_2$ with probabilities $1/2$, $%
1/4$, $1/4$. Further, the initial row of each matrix is zero, thus we may
consider only the lower-right $2\times 2$ submatrix: 
\[
\left( 
\begin{array}{cc}
0 & 0 \\ 
1 & 2
\end{array}
\right) =\left( 
\begin{array}{c}
0 \\ 
1
\end{array}
\right) \left( 
\begin{array}{cc}
1 & 2
\end{array}
\right) =\alpha _1\beta _1^T, 
\]
\[
\left( 
\begin{array}{cc}
4 & 0 \\ 
0 & 1
\end{array}
\right) =M, 
\]
\[
\left( 
\begin{array}{cc}
4 & 4 \\ 
0 & 0
\end{array}
\right) =\left( 
\begin{array}{c}
1 \\ 
0
\end{array}
\right) \left( 
\begin{array}{cc}
4 & 4
\end{array}
\right) =\alpha _2\beta _2^T. 
\]
Let $p_1=1/2$, $p_2=1/4$ denote the weights corresponding to $\alpha _1\beta
_1^T$, $\alpha _2\beta _2^T$ and $q=1/4$ denote the weight corresponding to $%
M$. The case of two $2\times 2$ matrices, one with rank $1$ and the other
with rank $2$, was solved in \cite{Pi, LR}. Our case involves three
matrices, two with rank $1$ and one with rank $2$, as well as the scaling
factor $1/2$ due to rewording. Generalizing, we obtain 
\[
\lambda =\frac 12\sum_{i=1}^2\sum_{j=1}^2\dsum\limits_{k=0}^\infty
p_ip_jq^k\ln \left( \beta _j^TM^k\alpha _i\right) . 
\]
Observe that 
\[
\beta _j^TM^k\alpha _i=\left\{ 
\begin{array}{ccc}
2 &  & \text{if }i=1\text{, }j=1 \\ 
2^2 &  & \text{if }i=1\text{, }j=2 \\ 
2^{2k} &  & \text{if }i=2\text{, }j=1 \\ 
2^{2(k+1)} &  & \text{if }i=2\text{, }j=2
\end{array}
\right. 
\]
and therefore 
\[
\frac{2\lambda }{\ln (2)}=\frac 14\dsum\limits_{k=0}^\infty \frac
1{4^k}+\frac 18\dsum\limits_{k=0}^\infty \frac 2{4^k}+\frac
18\dsum\limits_{k=0}^\infty \frac{2k}{4^k}+\frac
1{16}\dsum\limits_{k=0}^\infty \frac{2(k+1)}{4^k}=1, 
\]
as was to be shown. An independent proof of the conjecture was found by
Thomas Doumenc.

Let $\#(n)$ denote the number of $1$s in the binary expansion of $n$. We
know that $\#(n)=\ln (f(n))/\ln (2)$ and 
\[
\dsum\limits_{k=0}^n\#(k)\sim \frac 1{2\ln (2)}n\ln (n) 
\]
as $n\rightarrow \infty $. A corollary of our proof is that 
\[
\dsum\limits_{k=0}^{\left\lfloor n/3\right\rfloor }\#(3k)\sim \frac 1{2\ln
(2)}\frac n3\ln (n) 
\]
because $\#(3n)=\ln (g_3(n))/\ln (2)$. The formulas 
\[
\dsum\limits_{k=0}^{\left\lfloor n/3\right\rfloor }\#(3k+1)\sim \frac 1{2\ln
(2)}\frac n3\ln (n)\sim \dsum\limits_{k=0}^{\left\lfloor n/3\right\rfloor
}\#(3k+2) 
\]
follow similarly, that is, counting binary $1$s is (on average) independent
of ternary residue. We wonder whether simpler proofs of this fact can be
found.

\end{document}